\documentclass[11pt,oneside]{article}
\usepackage{amssymb,eucal}
\usepackage{amsfonts,amsmath,epsfig}

\usepackage[active]{srcltx}
\usepackage[sort&compress]{natbib}
\usepackage[colorlinks]{hyperref}

\numberwithin{equation}{section}
\newtheorem{theorem}{Theorem}[section]
\newtheorem{prop}[theorem]{Proposition}
\newtheorem{lemma}[theorem]{Lemma}
\newtheorem{assumption}{Condition}

\def\qed{\hfill$\Box$\medbreak}

\newenvironment{proof}[1][Proof]
{
  \medbreak\noindent\textit{#1.\/}
}{
  \qed
}

\newenvironment{proof*}[1][Proof]
{
  \medbreak\noindent\textit{#1.\/}
}{
}

\long\def\comment#1{{}}

\def\r{\epsilon}
\def\rr{\eta}
\def\gv{\,|\,}
\def\Grp#1{\left(#1\right)}
\def\Cbr#1{\left\{#1\right\}}
\def\Sbr#1{\left[#1\right]}

\def\Abs#1{\left|#1\right|}
\def\normx#1{\|#1\|}

\def\angx#1#2{\langle{#1},\,{#2}\rangle}

\def\Sp #1{\sp{({#1})}}
\def\nth#1{\frac{1}{#1}}

\def\tp{\sp{\top}}
\def\vf#1{\boldsymbol{#1}}
\def\eno#1#2{{#1}_1, \ldots, {#1}_{#2}}

\def\mean{\text{\sf E}}
\def\var{\text{\sf Var}}
\def\prob{\text{\sf Pr}}

\def\Reals{\mathbb{R}}
\def\Coms{\mathbb{C}}
\def\Ints{\mathbb{Z}}
\def\Nats{\mathbb{N}}

\def\Lsup{\mathop{\overline{\lim}}}

\def\toi{\to\infty}

\def\Lsup{\mathop{\overline{\mathrm{lim}}}}
\def\sppt{\mathrm{spt}}

\def\argx#1{\mathop{\arg#1}}

\def\rc{\varrho}  % Radius of convergence

\def\pset#1#2{\mathcal{D}\Grp{#1,#2}}  % parameter set
\def\grid{G}

\def\ball{B}
\def\parsp{{\Reals^p}}    % Space of parameters
\def\obssp{{\Reals^n}}    % Space of observations
\def\est#1{\widehat{#1}}
\def\btw#1#2#3{{#2\le#1\le#3}}

\def\diam{\delta}
\def\sfA{\text{\sf A}}
\def\sfC{\text{\sf C}}
\def\sfb{\text{\sf b}}
\def\sfd{\text{\sf d}}
\def\sfn{\text{\sf n}}
\def\sfr{\text{\sf r}}

\def\cD{\mathcal{D}}
\def\cE{\mathcal{E}}
\def\cJ{\mathcal{J}}
\def\cO{\mathcal{O}}
\def\cN{\mathcal{N}}
\def\cV{\mathcal{V}}

\begin{document}
\begin{center}
  \large
  \bf
  $\vf{L_0}$ regularized estimation for nonlinear models that have
  sparse underlying linear structures  \\[1.5ex]
  \rm\normalsize
  Zhiyi Chi
  \\
  Department of Statistics \\
  University of Connecticut \\
  215 Glenbrook Road, U-4120 \\
  Storrs, CT 06269, USA \\
  Email: zchi@stat.uconn.edu \\[1.5ex]
  \today
\end{center}

\begin{abstract}
  We study the estimation of $\beta$ for the nonlinear model $y =
  f(X\tp\beta) + \r$ when $f$ is a nonlinear transformation that is
  known, $\beta$ has sparse nonzero coordinates, and the number of
  observations can be much smaller than that of parameters ($n\ll p$).
  We show that in order to bound the $L_2$ error of the $L_0$
  regularized estimator $\est\beta$, i.e., $\normx{\est\beta -
    \beta}_2$, it is sufficient to establish two conditions.  Based on
  this, we obtain bounds of the $L_2$ error for (1) $L_0$ regularized
  maximum likelihood estimation (MLE) for exponential linear models
  and (2) $L_0$ regularized least square (LS) regression for the more
  general case where $f$ is analytic.  For the analytic case, we rely
  on power series expansion of $f$, which requires taking into account
  the singularities of $f$.

  \medbreak\noindent

  \medbreak\noindent

  \medbreak\noindent
  \textit{Acknowledgement.\/}  Research partially supported by
  NSF Grant DMS-07-06048 and NIH Grant MH-68028.

\end{abstract}
\section{Introduction}
Regularized estimation for sparse models that have a large number
of parameters comparing to that of observations has become an
important topic in statistics, machine learning, and a few other
areas \citep{efron:etal:04, zhao:yu:06, donoho:etal:06, bunea:etal:07,
  candes:tao:07, natarajan:95, field:94}.  The research in these areas
has been focused on regularized least square (LS) regression
for sparse linear models $y = X \beta + \r$, where $y\in\obssp$ is the
response vector, $X\in\Reals^{n\times p}$ the design matrix,
$\beta\in\parsp$ the vector of parameters, and $\r\in \obssp$ the
random error vector that has mean 0 given $X$.  By sparse we mean the
number of nonzero coordinates of $\beta$ is much smaller than $p$
\citep{wasserman:09}.

On the other hand, nonlinear models such as logistic models that have
underlying linear structures are widely used.  The general form of
such models is
\begin{align} \label{eq:nonlinear}
  y = f(X\tp \beta) + \r,
\end{align}
where $f: \Reals\to\Reals$ is a nonlinear function that may or may not
be known.  Here and henceforth, for $x = (\eno x n)\tp\in\obssp$, we
denote
\begin{align*}
  f(x) = (f(x_1), \ldots, f(x_n))\tp.
\end{align*}

The need for nonlinear models with sparse underlying linear structure
is clearly laid out in several recent works in neuroscience
\citep{sharpee:etal:08, sharpee:etal:04} and some algorithms based on
information criteria have been proposed to estimate not only $\beta$
but also $f$.  However, at this point, it seems very hard to evaluate
the estimation precision of those algorithms.

In this article we are content to establish the $L_2$ precision of
$L_0$ regularized estimator of $\beta$ for sparse models, when the
design matrix $X$ is fixed and $f$ is \emph{known\/}.  We shall
allow $n\ll p$.  Despite its limitation from a computational point of
view, the $L_0$ regularization is an important and conceptually simple
instrument for parameter estimation and model selection
\citep{akaike:74, schwarz:78, huang:etal:08}.  Besides, since many
improvements over the $L_0$ regularization are achieved by taking
advantage of properties of linear models that may fail to be had by
nonlinear models \citep{zhao:yu:06}, it is reasonable to take $L_0$
regularization as a prototype for further study on nonlinear models.
With this in mind, our concern is whether good estimation precision
\emph{could\/} be achieved instead of \emph{how fast\/} to achieve
it.

In Section \ref{sec:prelim}, we establish a basic result.  We show
that provided two conditions are satisfied, the $L_2$ error of the
$L_0$ regularized estimator satisfies a quadratic inequality which
yields the estimation precision.  Consequently, establishing the
estimation precision is reduced to establishing the two conditions.
As a minor benefit of the result, independence of the coordinates of
$\r$ in general need not be assumed.

We will also set up notation and collect other preliminary results in
Section \ref{sec:prelim}.  After that, we shall establish the alluded
conditions for exponential linear models and for analytic models,
i.e., models with analytic $f$.  Although a special case of analytic
models, exponential liner models are much simpler to handle due to its
explicit expression of the conditional density of $y$ given $X$.  For
these models, we consider the maximum likelihood estimator (MLE).  The
discussion is in Section \ref{sec:MLE}.  For analytic models, we will
consider the LS regression.  Sections \ref{sec:LSE-prelim} and
\ref{sec:LSE-A1} establish the two conditions, respectively.  In
Section \ref{sec:LSE-A1}, the approach is to use infinite power series
expansion of $f$.  The main complexity of the approach arises when $f$
has singularities on $\Coms$.  To illustrate, we will use as working
examples the logistic regression model in Section \ref{sec:MLE} and a
noise corrupted version of it in Section \ref{sec:LSE-A1}.  Most of
the proofs are collected in Section \ref{sec:proofs}.

\section{Preliminaries} \label{sec:prelim}
\subsection{Notation}
Denote by $X_1\tp$, \ldots, $X_n\tp$ the row vectors of $X$, with
$X_i \in\obssp$.  Denote by $\eno V p$ the column vectors of $X$.
We shall always assume that $X$ is fixed and impose the condition that
$V_j\not=0$.  In fact, if a column vector of $X$ is 0, then it has no
effect on $y$ and should be removed.  In the subsequent discussion,
the column vectors of $X$ should be understood as unnormalized.  It is
therefore helpful to think of $X$ as a collection of covariate
vectors registered exactly as they are observed.

For $S=\{i_1, \ldots, i_k\}$, with $1\le i_1 < \ldots < i_k\le p$,
denote $X_S=(V_{i_1}, \ldots, V_{i_k})$, and for $u \in\parsp$, denote
$u_S=(u_{i_1}, \ldots, u_{i_k})\tp$.  The support of $u$ is
\begin{align*}
  \sppt(u) =\{i: u_i\not=0\}.
\end{align*}
Denote by $\normx{u}_p$ the $L_p$ norm of $u$.  If $A$ is a set, denote
by $|A|$ its cardinality.  The $L_0$ norm of $u$ refers to
$|\sppt(u)|$ and is often denoted by $\normx{u}_0$.  We choose the
notation $|\sppt(u)|$ since it seems more intuitive.

For $\varphi=(\varphi_1, \ldots, \varphi_n)$ and $x\in
\obssp$, where each $\varphi_i: \Reals\to \Reals$, denote
\begin{align*}
  \varphi(x) = (\varphi_1(x_1), \ldots, \varphi_n(x_n))\tp.
\end{align*}

\subsection{General form of estimator and line of argument}  
The general form of an $L_0$ regularized estimator is
\begin{align} \label{eq:general}
  \est\beta = \argx\min_{u\in D} \Sbr{\ell(y, Xu) + c_r|\sppt(u)|},
\end{align}
where $D$ is a pre-selected search domain in $\parsp$, $\ell(y, Xu)$
is certain loss function, and $c_r>0$ is a tuning
parameter.  For the MLE, $\ell(y, Xu)$ is the minus log likelihood,
while for the LS regression, it is $\normx{y-Xu}_2^2$.  For linear
regression, $D$ is typically set equal to $\parsp$.  However, for
nonlinear regression, our position is that some constraint on $D$ is
needed in order to control the potentially large variation of the
functional property of $f$ at different possible values of
$X\beta$.

For both the MLE and LS regression, the argument to establish
the precision of $\est\beta$ proceeds as follows.  First,
it is easy to show that $\est\beta$ satisfies an inequality
of the following 
form,
\begin{align} \label{eq:ineq}
  G(\psi(X\est\beta) - \psi(X\beta))
  \le 2 |\angx{\r}{\varphi(X\est\beta) - \varphi(X\beta)}| 
  -
  c_r(|\sppt(\est\beta)|-|\sppt(\beta)|),
\end{align}
where $G$ is a function $\obssp\to\Reals$, $\psi = (\eno\psi n)$ and
$\varphi=(\eno \varphi n)$, with $\psi_i$ and $\varphi_i$ being
functions $\Reals\to\Reals$.  Then the following two conditions will
be established.

\renewcommand{\theassumption}{H\arabic{assumption}}
\begin{assumption}\rm\label{cond:A1}
  Given $q\in (0,1)$, there is $c_1=c_1(X, \beta, \varphi,q)>0$, such
  that
  \begin{align*}
    \prob\Cbr{
      |\angx{\r}{\varphi(Xu) - \varphi(X\beta)}|
      \le c_1\sqrt{n} \normx{u-\beta}_1, 
      \text{ all } u \in D
    }\ge 1-2q.
  \end{align*}
  The coefficient 2 in $1-2q$ is nonessential.  It is for ease of
  notation in the statements of main results.
\end{assumption}

\begin{assumption}\rm \label{cond:A2}
  There is $c_2=c_2(X, \beta, \psi)>0$, such that for all $u\in D$,
  \begin{align*}
    G(\psi(X u) - \psi(X\beta)) \ge c_2 n \normx{u-\beta}_2^2.
  \end{align*}
\end{assumption}

The constants $c_1$ and $c_2$ will be explictly constructed.  In
general, both depend on $X$.  Since we only consider fixed design,
they are nonrandom.

We will check the conditions respectively for the MLE and LS
regression.  Once this is done, using the next result, we then obtain
a bound on $\normx{\est\beta - \beta}_2$.  Note that the result is
stated in a little more general form as it does not require that
$\est\beta$ be the one defined by \eqref{eq:general}.

\begin{prop} \label{prop:basic}
  Suppose \emph{Conditions \ref{cond:A1} and \ref{cond:A2}\/} are
  satisfied.  If $\est\beta\in D$ is a random variable that always
  satisfies the inequality \eqref{eq:ineq} with $c_r = 3 c_1^2/c_2$,
  then, letting $\kappa_r = 3 c_1/c_2$,
  \begin{align*}
    \prob\Cbr{\normx{\est\beta-\beta}_2 \le
      \frac{\kappa_r\sqrt{|\sppt(\beta)|}}{\sqrt{n}}
    } \ge 1-2q.
  \end{align*}
\end{prop}

In order for the bounds to be meaningful, we need to make sure
$\kappa_r$ is not too large, at least comparing to $\sqrt{n}$.  This
will be the main consideration when we try to establish Conditions
\ref{cond:A1} and \ref{cond:A2}.

Because Proposition \ref{prop:basic} plays a fundamental role in our
study, we give its proof below.  This is the only result whose proof
appears in the main text.

\begin{proof}[Proof of Proposition \emph{\ref{prop:basic}\/}]
  Denote $T=\sppt(\beta)$ and $S=\sppt(\est\beta)$.  Under Conditions
  \ref{cond:A1} and \ref{cond:A2}, with probability at least $1-2q$, 
  \begin{align*}
    c_2 n \normx{\est\beta - \beta}_2^2
    &
    \le 2 c_1 \sqrt{n} \normx{\est\beta - \beta}_1
    - c_r (|S|-|T|) \\
    &
    \le 2 c_1 \sqrt{n} \sqrt{|S\cup T|} \normx{\est\beta - \beta}_2
    - c_r (|S|-|T|),
  \end{align*}
  where the second inequality is due to $\sppt(\beta-\est\beta)
  \subset S\cup T$ and Cauchy-Schwartz inequality.  Let 
  $t = \normx{\est\beta - \beta}_2$ and $b=c_1/c_2$.  Then
  \begin{align*}
    t^2 
    - \frac{2 b \sqrt{|S\cup T|} t}{\sqrt{n}}+
    \frac{3 b^2(|S|-|T|)}{n}\le
    0.
  \end{align*}
  The left hand side is a quadratic function in $t$.  In order for the
  inequality to hold, there have to be $|S\cup T|\ge 3(|S|-|T|)$ and
  \begin{align*}
    0\le t \le \frac{b}{\sqrt{n}} 
    \Sbr{\sqrt{|S\cup T|} + \sqrt{|S\cup T| + 3(|T|-|S|)}\,}.
  \end{align*}

  Let $T_1 = T\setminus S$ and $S_1 = S\setminus T$.  By $|S\cup T| =
  |S_1|+|T|$ and $|T|-|S|=|T_1|-|S_1|$,
  \begin{align*}
    0\le t \le \frac{b}{\sqrt{n}} 
    \Grp{\sqrt{|T|+|S_1|} + \sqrt{|T| + 3|T_1|-2|S_1|}}.
  \end{align*}
  It is easy to see that due to $|T_1|\le |T|$, the right hand side
  is a decreasing function in $|S_1|$ on $[0, (|T|+3|T_1|)/2]$, and
  hence is no greater than its value at 0, which is
  $(b/\sqrt{n})(\sqrt{|T|} + \sqrt{|T|+3|T_1|}) \le 3b\sqrt{|T|} /
  \sqrt{n}$.
\end{proof}

To establish Conditions \ref{cond:A1} and \ref{cond:A2}, certain
assumptions are needed.  We next discuss the major assumptions used by
both the MLE and LS regression.

\subsection{Tail assumption on errors}
To establish Condition \ref{cond:A1}, we will need the following
assumption on $\r$. 

\paragraph{Tail assumption.}  There is $\sigma>0$, such that for any
$t$, $\eno a n$ $\in$ $\Reals$,
\begin{align} \label{eq:tail}
  \prob\Cbr{
    \Grp{\sum_{i=1}^n a_i\r_i}^2>t^2\sum_{i=1}^n a_i^2
  } \le 2\exp\Cbr{-\frac{t^2}{2\sigma^2}}.
\end{align}

The tail assumption \eqref{eq:tail} rather mild.  If $\r\sim
N(0,\sigma^2\Sigma)$ and the spectral radius of $\Sigma$ is no greater 
than 1, then \eqref{eq:tail} holds.  In this case, $\eno\r n$ need not
be independent.  Moreover, if $\r_i$ are independent, such that
$\mean(\r_i)=0$ and $|\r_i|\le\sigma$ for all $i$, then by Hoeffding's
inequality \citep{pollard:84}, \eqref{eq:tail} holds.

\subsection{Coherence and restricted domains} \label{ssec:coherence}
In order to identify $\beta$, some conditions on the correlations
between the column vectors of $X$ are needed.  The maximum correlation
between columns of $X$ is
\begin{align*}
  \mu(X) = \sup_{1\le i<j\le p} \frac{|V_i\tp V_j|}{\normx{V_i}_2
    \normx{V_j}_2}.
\end{align*}
Conditions on $\mu(X)$ are often referred to as coherence property
\citep{bunea:etal:07, candes:plan:09}.  The following function
\begin{align} \label{eq:coherence}
  \sfn(\nu)=(1-\nu) \Sbr{1+1/\mu(X)}
\end{align}
will be regularly used in our discussion.
\begin{prop} \label{prop:separable}
  Fix $\nu\in [0,1]$.  (1) For $u\in\parsp$, if $|\sppt(u)|\le
  \sfn(\nu)$, then
  \begin{align*}
    \normx{Xu}_2^2 \ge \nu [1+\mu(X)]\sum_{j=1}^p |u_j|^2
    \normx{V_j}_2^2.
  \end{align*}
  (2) For $u$, $v\in\parsp$, if $|\sppt(u)\cup \sppt(v)|\le
  \sfn(\nu)$, then
  \begin{align*}
    \normx{X(u-v)}_2^2 \ge \nu [1+\mu(X)]\sum_{j=1}^p |u_j-v_j|^2
    \normx{V_j}_2^2.
  \end{align*}
  In particular, the inequality holds if $|\sppt(u)|\vee|\sppt(v)|\le
  \sfn(\nu)/2$.
\end{prop}

As mentioned earlier, for the estimator \eqref{eq:general}, we
need to impose some constraints on the search domain $D$.  For this
purpose, we define several sets.  For $I\subset\Reals$, let
\begin{align}
  \cD(I)
  = 
  \{u\in\parsp: X_i\tp u\in I, \ \btw i 1 n\}, \label{eq:domain}
\end{align}
and for $h\ge 1$, let
\begin{align}
  \pset I h
  =\cD(I) \cap \Cbr{u\in \parsp: |\sppt(u)|\le h}.
  \label{eq:domain2}
\end{align}
Apparently, denoting by $T$ the mapping $u\to Xu$, $\cD(I)=
T^{-1}(I^n)$.

One constraint that will be regularly imposed is $D\subset \pset I
{\sfn(\nu)/2}$ for some $\nu\in (0,1)$.  The implied constraint that
$X_i\tp u\in I$ for every $i$ is to make sure that the functions
involved in the estimator \eqref{eq:general}, i.e., $G$, $\psi_i$ and
$\phi_i$, have good enough properties for all candidate values of
$\beta$, especially properties determined by derivatives.  This
constraint on the functional properties is needed when we
establish both Conditions \ref{cond:A1} and \ref{cond:A2}.  For linear
regression, roughly speaking, this is not a concern and one can simply
choose $I=\Reals$, simply because the derivative of a linear function
is constant, and so the pertinent functional properties are uniform.

The constraint $D\subset \pset I{\sfn(\nu)/2}$ also imposes a
constraint on $|\sppt(\est\beta)|$.  As Proposition
\ref{prop:separable} indicates, one consequence of the constraint is
that any two candidate estimates of $\beta$ can be well separated by
their corresponding values of $Xu$, so that a large portion of $\beta$
can be correctly identified.  For this reason, the constraint will be
needed when we establish Condition \ref{cond:A2}.  Clearly, the
smaller $\mu(X)$ is, the 
milder the constraint.  Under mild conditions, $\mu(X)$ can
be as small as $O(\sqrt{n^{-1}\ln p})$; see \cite{candes:plan:09} and
also the comments at the end of Section \ref{ssec:log-r}.  This
results in a constraint of the form $|\sppt(\est\beta)|\le
C\sqrt{n/\ln p}$, which is quite mild even when $p$ is much larger
than $n$, for example, $p=n^a$ for some $a>1$.

We shall need the following properties of $\pset I h$.
\begin{prop} \label{prop:compact}
  (1) If $I$ is closed, then $\pset I 1\subset \pset I 2 \subset
  \cdots$ are closed and (2)~if $I$ is compact and $h < \sfn(0)=
  1+\mu(X)^{-1}$, then $\pset I h$ is compact.
\end{prop}

\section{Exponential linear models}  \label{sec:MLE}
\subsection{Setup and main result} \label{ssec:glm-setup}
Let $\mu$ be a Borel measure on $\Reals$ with $\mu(\Reals)>0$.
Suppose $I\subset\Reals$ is an nonempty open interval and $\{P_t:
t\in I\}$ is a family of probability distributions on $\Reals$, such
that with respect to $\mu$ each $P_t$ has a density
\begin{align} \label{eq:glm-form}
  p_t(y) = \exp\Cbr{t y - \Lambda(t)},
  \ \text{with} \ 
  \Lambda(t) = \ln \Sbr{\int e^{t y}\,\mu(dy)}.
\end{align}
As is well known, $\Lambda\in C^\infty(I)$ and for $t\in I$,
\begin{align} \label{eq:glm-moments}
  \mean(\xi) = \Lambda'(t), \quad \var(\xi) =
  \Lambda''(t)>0, \quad\text{if } \xi\sim P_t.
\end{align}
For example, if $\mu=N(0,\sigma^2)$, then $\Lambda(t) = \sigma^2
t^2/2$ and $P_t=N(\sigma^2 t, \sigma^2)$.  If $\mu$ is the counting
measure on $\{0,1\}$, then $\Lambda(t) = \ln(1+e^t)$ and $P_t$ is the
Bernoulli distribution with parameter $e^t/(1+e^t)$.  We notice that
given $y$, $g(t):=p_t(y)$ can be ananlyticall extended to the domain
$\{z\in\Coms:{\rm Re}(z)\in I\}$.  This fact is not needed in the
rest of the section.

Assume that given $X$, $\eno y n$ are independent, such that each
$y_i\sim P_{t_i}$ with $t_i = X_i\tp\beta$.  The joint likelihood of
$\eno y n$ is then
\begin{align*}
  \prod_{i=1}^n \exp\Cbr{y_i X_i\tp\beta -
    \Lambda(X_i\tp\beta)}
  =
  \exp\Cbr{y\tp X\beta - \sum_{i=1}^n \Lambda(X_i\tp\beta)}.
\end{align*}
From the expression, the $L_0$ regularized MLE for $\beta$ is 
\begin{align} \label{eq:MLE}
  \est\beta = \argx\max_{u\in D}
  \Sbr{
    y\tp Xu - \sum_{i=1}^n \Lambda(X_i\tp u) - c_r|\sppt(u)|
  }.
\end{align}
If $\beta\in D$, then  
\begin{align*}
  y\tp X\beta - \sum_{i=1}^n \Lambda(X_i\tp\beta) - 
  c_r|\sppt(\beta)|
  \le
  y\tp X\est\beta - \sum_{i=1}^n \Lambda(X_i\tp\est\beta) - 
  c_r|\sppt(\est\beta)|,
\end{align*}
and hence
\begin{align*}
  &
  \sum_{i=1}^n\Sbr{
    \Lambda(X_i\tp\est\beta) - \Lambda(X_i\tp\beta) -
    \Lambda'(X_i\tp\beta)X_i\tp(\est\beta-\beta)
  } \\
  &
  \le \angx{\r}{X\est\beta- X\beta} - c_r(|\sppt(\est\beta)| -
  |\sppt(\beta)|),
\end{align*}
where $\r_i = y_i - \mean(y_i) = y_i - \Lambda'(X_i\tp\beta)$ has mean
0 for each $i$.  It is seen that the inequality gives rise to
\eqref{eq:ineq} once we
define
\begin{align}
  G(x)=\sum_{i=1}^n x_i,\quad
  \psi_i(z)
  = \Lambda(z) - \Lambda'(X_i\tp\beta) z, \quad
  \varphi_i(z) =z/2,
  \label{eq:glm-f}
\end{align}
for $x\in\obssp$, $z\in\Reals$ and $\btw i 1 n$.

\begin{theorem} \label{thm:glm}
  Suppose $\eno\r n$ satisfy \eqref{eq:tail} for some $\sigma>0$.
  Fix $\nu\in(0,1)$.  Let $D =\pset I {\sfn(\nu)/2}$ in
  \eqref{eq:MLE}, where $\sfn(\nu)$ is defined in
  \eqref{eq:coherence}.  Suppose
  \begin{align} \label{eq:glm-db}
    \delta:= \inf_{t\in I} \Lambda''(t)>0.
  \end{align}
  Fix $q\in (0,1/2)$.  Let
  \begin{align*}
    c_r = \frac{3\sigma^2 \ln(p/q)}{\nu\delta[1+\mu(X)]} 
    \frac{\max_j \normx{V_j}_2^2}{\min_j \normx{V_j}_2^2}
  \end{align*}
  in \eqref{eq:MLE}.  Then, provided $\beta\in D$,
  \begin{align} \label{eq:glm-est}
    &\prob\Cbr{
      \normx{\est\beta-\beta}_2
      \le \frac{\kappa_r\sqrt{|\sppt(\beta)|}}{\sqrt{n}}
    } \ge 1-2q, \\
    \text{where}\quad
    &
    \kappa_r = \frac{3\sigma\sqrt{2\ln(p/q)}}{\nu\delta[1+\mu(X)]}
    \times \frac{\sqrt{n}\max_j \normx{V_j}_2}{
      \min_j\normx{V_j}_2^2}.
    \nonumber
  \end{align}
\end{theorem}

\subsection{Comments} \label{ssec:glm-notes}
Some comments on Theorem \ref{thm:glm} are in order, many of them
also apply to the results we shall establish later.  First, on the
constraint $\est\beta\in \pset I {\sfn(\nu)/2}$.  As noted in Section
\ref{ssec:coherence}, under mild conditions, for $p$ with $\ln p =
o(n)$, $\sfn(\nu)\asymp\sqrt{n/\ln   p}$.  In many cases, since it is
reasonable to assume that $|\sppt(\beta)|=O(1)$ \citep{wasserman:09},
the constraint then is very mild.

Second, on $\normx{\est\beta-\beta}_2$, which is determined by
$\kappa_r\sqrt{|\sppt(\beta)}/\sqrt{n}$ in \eqref{eq:glm-est}.  By
\eqref{eq:glm-est}, $\kappa_r=O(R\sqrt{\ln p})$, where
\begin{align*}
  R
  =\frac{\sqrt{n} \max_j \normx{V_j}_2}{\min_j \normx{V_j}_2^2}
  =\frac{\max_j \normx{V_j}_2/\sqrt{n}}{\min_j \normx{V_j}_2^2/n}.
\end{align*}
Under mild conditions, $R$ grows very slowly with $n$.  For example,
$R=1$ if $X$ is such that $\normx{V_j}_2 = \sqrt{n}$ (recall all
$V_j\in\obssp$).  We shall see such an example related to the logistic
regression.  As another example, suppose all the $np$ entries of $X$
are i.i.d. $\sim Z$.  If $Z$ is bounded, then clearly $\max_j
\normx{V_j}_2/\sqrt{n} = O(1)$.  If $Z\sim N(0,1)$, then for any
$0<\rr<1/2$,
\begin{align*}
  \prob\Cbr{\max_{\btw j 1 p} \normx{V_j}_\infty \le \sqrt{2\ln
      (np/\rr)}} \ge 1-2\rr.
\end{align*}
Since $\max_j\normx{V_j}_2 \le \sqrt{n} \max_j\normx{V_j}_\infty$,
then with high probability, $\max_j \normx{V_j}_2/\sqrt{n}
 = O(\sqrt{\ln (np)})$.  At the same time, given $0<c< \mean(Z^2)$, 
\begin{align*}
  \prob\Cbr{\nth n\min_{\btw j 1 p} \normx{V_j}_2^2 \le c}
  \le p \prob\Cbr{Z_1^2+\cdots + Z_n^2 \le n c}
  \le p \psi(c)^n,
\end{align*}
where $\psi(c) = \inf_{t>0} \mean[e^{t c - tZ^2}] < 1$.  Therefore,
for large $n$ and $p$, with high
probability, we have $\max_j \normx{V_j}_2/\sqrt{n} =
O(\sqrt{\ln(np)})$ or even $O(1)$ on the one hand, and $\min_j
\normx{V_j}_2^2/n \ge c$ on the other, provided $\ln p = o(n)$.  In
particular, suppose $p=O(n^a)$ for some $a>0$.  Then it is seen that
$R = O(\sqrt{\ln n})$ or even $O(1)$, and hence, by
\eqref{eq:glm-est}, with high probability,
$\normx{\est\beta-\beta}_2=O(\ln n/\sqrt{n})$ or $O(\sqrt{\ln
  p}/\sqrt{n})$.

Finally, the precision also depends on $\delta=\inf_{t\in I}
\Lambda''(t)$.  To see why $\delta$ matters, consider the case where
$\Lambda''(t)$ is uniformly small in an interval $I$ that contains all
of $X_i\tp \beta$.  This implies that $\Lambda'(t)$ has little change
on $I$, so by \eqref{eq:glm-moments}, $\mean(y_1)$, \ldots,
$\mean(y_n)$ are close to each other, and at the same time each $y_i$
has little variation.  This gives rise to a nearly ``flat'' plot of
$y_i$ \emph{vs\/} $X_i\tp\beta$, which makes the 
identification of $\beta$ difficult.  That is to say the precision
of the estimate cannot be high.  Certainly, if $\Lambda''(t)$ has a
wide range on $I$, then using $\inf_{t\in I} \Lambda''(t)$ to set
$c_r$ can be quite conservative.  However, as $X_i\tp\beta$ are
unknown, it is the only way to account for all the possible values of
$X_i\tp\beta$, including the least ideal one.

\subsection{Logistic regression} \label{ssec:log-r}
Suppose $\eno y n$ are independent Bernoulli random variables, such
that
\begin{align*}
  \prob\Cbr{y_i=1} 
  = e^{X_i\tp\beta}/(1+e^{X_i\tp \beta}), \quad
  i=1,\ldots,n.
\end{align*}
The corresponding parametric family of densities is $p_t(y) = \exp\{t
y - \Lambda(t)\}$ with respect to the counting measure on $\{0,1\}$,
with $\Lambda(t) = \ln(1+e^t)$.

For $i=1,\ldots,n$, $\r_i = y_i - \prob\Cbr{y_i=1} \in (-1,1)$.
Therefore, by Hoeffding's inequality \citep{pollard:84},
\eqref{eq:tail} holds with $\sigma=1$.  Given $I\subset \Reals$, by
direct calculation,
\begin{align*}
  \inf_{t\in I} \Lambda''(t)
  = \Grp{2\cosh\frac{M_I}{2}}^{-2},
  \ \text{with}\ M_I =\sup_{t\in I} |t|.
\end{align*}
Given $q\in (0,1)$, let
\begin{align*}
  c_r =  \frac{12 \ln (p/q)}
  {\nu[1+\mu(X)]}
  \times \frac{\max_j \normx{V_j}_2^2}{\min_j \normx{V_j}_2^2}
  \times \cosh^2\frac{M_I}{2}
\end{align*}
and 
\begin{align*}
  \kappa_r = \frac{12\sqrt{2\ln(p/q)}}{\nu[1+\mu(X)]}
  \times \frac{\sqrt{n}\max_j\normx{V_j}_2}
  {\min_j\normx{V_j}_2^2} \times \cosh^2\frac{M_I}{2}.
\end{align*}
By Theorem \ref{thm:glm}, if $\beta\in \pset I {\sfn(\nu)/2}$, then,
with probability at least $q$, \eqref{eq:glm-est} holds for the
estimator
\begin{align*}
  \est\beta=\argx\max\Cbr{
    y\tp X u - \sum_{i=1}^n \ln (1+e^{X_i\tp u}) - 
    c_r|\sppt(u)|:  u\in \pset I {\sfn(\nu)/2}
  }.
\end{align*}

If $X$ is binary, i.e., $X_{ij}=0$ or $1$, the result can be somewhat
simplified.  Let $\tilde X\in \Reals^{n\times (p+1)}$
such that $\tilde X_{ij} = 2 X_{ij}-1$, for $j\le p$ and $\tilde
X_{i,p+1}=1$.  Also let $\tilde\beta\in\Reals^{p+1}$ such that
$\tilde\beta_j = \beta_j/2$ for $j\le p$ and $\tilde\beta_{p+1} =
\sum_{j=1}^p \beta_j/2$.  Then $X_i\tp \beta = \tilde X_i\tp
\tilde\beta$.  Let $\eno {\tilde V}{p+1}$ be the column vectors of
$\tilde X$.  Then $\normx{\tilde V_j}_2=\sqrt{n}$.  If we regress $y$
on $\tilde X$ to estimate $\tilde\beta$, then
\begin{align*}
  c_r = \frac{12\ln [(p+1)/q]}{\nu[1+\mu(\tilde X)]} \times
  \cosh^2\frac{M_I}{2}, \quad
  \kappa_r = \frac{12\sqrt{2\ln(p/q)}}{\nu[1+\mu(\tilde X)]}
  \times \cosh^2\frac{M_I}{2}.
\end{align*}

In the example, $\mu(\tilde X)$ can be very small.  If $X_{ij}$ are
i.i.d.\ with $\prob\{X_{ij}=0\}=\prob\{X_{ij}=1\}=1/2$, then for
\emph{any\/} $1\le j<k\le p+1$, $\tilde V_j\tp \tilde V_k\sim
\sum_{i=1}^n \rr_i$, where $\rr_i$ are i.i.d.\ with $\prob\{\rr_i=1\}
= \prob\{\rr_i=-1\}=1/2$.  By Hoeffing's inequality, given $t>0$,
\begin{align*}
  \prob\Cbr{\frac{|\tilde V_j\tp \tilde V_k|}
    {\normx{\tilde V_j}_2 \normx{\tilde V_k}_2}
    \ge \frac{t}{\sqrt{n}}
  }
  =
  \prob\Cbr{\Abs{\sum_{i=1}^n \rr_i}\ge t\sqrt{n}}
  \le 2 e^{-t^2/2}.
\end{align*}

It follows that given $\delta\in (0,1)$, 
\begin{align*}
  &
  \prob\Cbr{\mu(\tilde X)\ge
    \sqrt{\frac{2}{n} \ln \frac{(p+1)^2}{\delta}}
  }\\
  &\le
  \frac{p(p+1)}{2}
  \prob\Cbr{\frac{|\tilde V_1\tp \tilde V_2|}
    {\normx{\tilde V_1}_2 \normx{\tilde V_2}_2}
    \ge 
    \sqrt{\frac{2}{n} \ln \frac{(p+1)^2}{\delta}}
  }
  \le \delta.
\end{align*}
Therefore, with high probability, $\mu(\tilde X)=O(\sqrt{\ln p/n})$,
which is very small for reasonably large $p$ and $n$.

\section{Least square regression:  preliminaries}
\label{sec:LSE-prelim}
\subsection{Reformulation and Condition \ref{cond:A2}}
Suppose that, with $X$ fixed,
\begin{align*}
  y_i = f(X_i\tp\beta)+\r_i, \quad \btw i 1 n,
\end{align*}
where $\r_i$ are independent with mean 0.  The $L_0$ regularized
LS estimator for $\beta$ is
\begin{align} \label{eq:LSE}
  \est\beta = \argx\min_{u\in D}
  \Sbr{
    \normx{y-f(Xu)}_2^2 + c_r|\sppt(u)|
  },
\end{align}
where, as in \eqref{eq:MLE}, $D$ is a suitable search domain in
$\parsp$ and $c_r$ is a regularization parameter.  If $\beta\in D$,
then
\begin{align*}
  \normx{y-f(X\est\beta)}_2^2 + c_r|\sppt(\est\beta)|\le
  \normx{y-f(X\beta)}_2^2 + c_r|\sppt(\beta)|,
\end{align*}
and hence
\begin{align*}
  \normx{f(X\est\beta) - f(X\beta)}_2^2 \le
  2\angx{\r}{f(X\est\beta) - f(X\beta)}-c_r(|\sppt(\est\beta)|
  - |\sppt(\beta)|),
\end{align*}
which implies \eqref{eq:ineq} once we define
\begin{align}
  G(x) = \normx{x}_2^2, \quad
  \psi_i(z) = \varphi_i(z) = f(z),
  \label{eq:lse-f}
\end{align}
for $x\in\obssp$, $z\in\Reals$ and $\btw i 1 n$.  By Proposition
\ref{prop:basic}, all we need to do then is to find suitable constants
$c_1$ and $c_2$ so that Conditions \ref{cond:A1} and \ref{cond:A2} are
satisfied.

For $I\subset \Reals$ that contains at least two points, denote
\begin{align*}
  \sfd(f,I)=
  \inf \Cbr{\frac{|f(x)-f(y)|}{|x-y|}:\, x\in I,\, y\in I,\,
    x\not=y}.
\end{align*}

We start with the easier task of establishing Condition \ref{cond:A2}.
\begin{prop} \label{prop:lse-d}
  Let $I\subset\Reals$ be an interval with positive length.  Suppose
  $f$ is defined on $I$ with $\sfd(f,I)>0$.  Fix $\nu\in (0,1)$.  Let
  $D$ in \eqref{eq:LSE} be a subset of $\pset I {\sfn(\nu)/2}$.  If
  $\beta\in D$, then for $G$ and $\psi$ defined as in
  \eqref{eq:lse-f}, \emph{Condition \ref{cond:A2}\/} is satisfied with
  \begin{align*}
    c_2 = \frac{\sfd(f,I)^2 \nu[1+\mu(X)]}{n}
    \min_{\btw j 1 p}\normx{V_j}_2^2.
  \end{align*}
\end{prop}

As noted in Section \ref{ssec:glm-notes}, under mild conditions, for
large $n$ and reasonably large $p$, $c_2\asymp 1$.  Therefore, by
Proposition \ref{prop:basic}, in order for the estimate $\est\beta$ to
have some reasonable precision, the coefficient $c_1$ in Condition
\ref{cond:A1} has to be of order $o(\sqrt{n})$.  To this end,
depending on how well the nonlinear function $f$ behaves, some extra
constraints need to be imposed on the domain $D$.  Section
\ref{sec:LSE-A1} is devoted to establishing Condition \ref{cond:A1}
for the LS regression.  Below we outline the steps to be
taken.

\subsection{Observations that point to Condition \ref{cond:A1}}
Recall that Condition \ref{cond:A1} stipulates an upper bound on
$|\angx{\r}{f(Xu)-f(X\beta)}|$ that has to hold
\emph{simultaneously\/} for all $u$.  If $f(x)=x$, such a bound is
easy to find due to the conjugate relation $\angx{\r}{f(Xu)-f(X\beta)}
= \angx{X\tp\r}{u-\beta}$, as it then suffices to find a bound for
$\normx{X\tp\r}_\infty$, which can be derived from the tail assumption
on $\r$ \citep{zhang:09, candes:plan:09}.  For nonlinear $f$, in
general, there are no similar applicable relations.  However, like
$e^x/(1+e^x)$, in many cases, $f$ is analytic and so we may exploit
its power series expansions around different points.  By working with,
say $f(x)=x^2$, one could imagine a kind of power series expansion
\begin{align*}
  f(Xu)=\sum M_\alpha h_\alpha(u),
\end{align*}
such that each $M_\alpha $ is some type of (row-wise) monomial
transformation of $X$, and $h_\alpha(u)$ a vector resulting from a
similar transformation of $u$.  This makes it possible to rewrite
$\angx{\r}{f(Xu) - f(X\beta)}$ as an infinite sum of
$\angx{M_\alpha\tp\r}{h_\alpha(u) - h_\alpha (\beta)}$, which could
lead to a desirable bound.

The method works if $f$ is analytic on the entire $\Coms$, or, more
generally, when all the coordinates of $X u$ and $X\beta$ fall into
the disc of convergence of the power series expansion of $f$ at 0. 
On the other hand, when $f$ has poles as $e^x/(1+e^x)$ does, the
coordinates of $Xu$ and $X\beta$ may fall into different discs of
convergence of power series expansion.   Roughly, to deal with this
problem, our approach is to cover the line segment connecting $Xu$ and
$X\beta$ with different discs of convergence of power series, apply
the result obtained for the case of single analytic disc, and patch
together the resulting bounds.  This turns out to account for most of
the complexity in our treatment of the analytic case.

One question is whether we can just use a finite Taylor
expansion to derive bounds for $\angx{\r}{f(Xu)-f(X\beta)}$, thus
dispensing with the assumption of analyticity.  The answer seems to be
no in general.  Unless $f$ is a polynomial, a finite Taylor expansion
of $f(Xu)-f(X\beta)$ has a remainder term of the form $R_\alpha(u)
[h_\alpha(u) - h_\alpha(\beta)]$, where $R_\alpha(u)$ is a matrix that
in general depends on $u$.  As a result, although for each individual
$u$, we can get a bound for $\angx{\r\tp R_\alpha(u)} {h_\alpha(u) -
  h_\alpha(\beta)}$ that holds with high probability, there is no
guarantee to get that with high probability, the bounds hold
simultaneously for all $u$, which is needed for establishing the
precision of $\est\beta$.

\section{Least square regression: continued}
\label{sec:LSE-A1}

\subsection{Setup}
Let $I\subset\Reals$ be a closed interval with positive length.  In
this section, we assume that $f: I\to\Reals$ is analytic in a
neighborhood of $I$, i.e., $f$ has a (unique) analytic extension onto
an open set in $\Coms$ containing $I$.  This is equivalent to
saying that $f \in C^\infty(I)$ and for each $t\in I$, there is $r>0$,
such that
\begin{multline} \label{eq:analytic}
  \sum_{k=0}^\infty |a_k| r^k < \infty,
  \ \text{where}\ \
  a_k = \frac{f\Sp k(t)}{k!}\in\Reals,\\
  \text{and}\quad
  f(z+t) = \sum_{k=0}^\infty a_k z^k, \ \text{for all} \
  z\in (-r, r) \text{ with } z+t\in I.
\end{multline}
The radius of convergence of the power series
\eqref{eq:analytic}, henceforth denoted by $\rc(f,t)$, can be
determined by \citep{rudin:87}
\begin{align*}
  \rc(f,t) = \Grp{\Lsup_{k\toi} |a_k|^{1/k}}^{-1}
\end{align*}
If $|z|<\rc(f,t)$, then we say $f(z+t)$ has a convergent power series
expansion at $t$.

We will regularly use the following weighted $L_1$ norm
\begin{align} \label{eq:w-L1-norm}
  \normx{u}_{1,s} = \sum_{j=1}^p |u_j|\normx{V_j}_s,
  \quad u\in\parsp, \ s\ge 1.
\end{align}
Recall that it is assumed from the beginning that $V_j\not=0$ for all
$j$.  Therefore, $\normx{u}_{1,s}$ is indeed a norm.  Finally, if
$(\cE, \normx{\cdot})$ is a normed linear space, then denote by
\begin{align*}
  \ball(u,a; \normx{\cdot}) = \{v\in \cE: \normx{v-u}<a\}
\end{align*}
the sphere centered at $u\in\cE$ with radius $a>0$ under the norm
$\normx{\cdot}$, and by 
\begin{align*}
  \diam(E; \normx{\cdot})
  = \inf\{a: E\subset \ball(u, a; \normx{\cdot}) \text{ for some } u\}.
\end{align*}
the infimum of the radii of spheres under the norm $\normx{\cdot}$
that contain $E\subset\cE$.

\subsection{Single analytic disc} \label{ssec:one-disc}

We first consider the case where all $f(X_1\tp u)$, \ldots, $f(X_n\tp
u)$ have convergent power series expansions at 0.  The main result of
this section is as follows.
\begin{theorem} \label{thm:one-disc}
  Suppose $0\in I$ and  $\sfd(f,I)>0$.  Fix $\nu\in (0,1)$ and
  $\theta\in (0,1)$.  Suppose
  \begin{align*}
    \begin{split}
      D=\pset I {\sfn(\nu)/2} \cap 
      \Cbr{u\in\parsp: \normx{u}_{1,\infty}\le \theta\rc(f,0)/2}
    \end{split}
  \end{align*}
  in \eqref{eq:LSE} and $\r$ satisfies \eqref{eq:tail} for $\sigma>0$.
  Given $q\in (0,1)$, let $\lambda_p = \ln [p(1+q^{-1})]$.   If
  $\beta\in D$, then the conclusion of \emph{Proposition
    \ref{prop:basic}\/} holds with
  \begin{align*}
    c_1 = \sigma\sqrt{2\lambda_p}
    \sum_{k=1}^\infty
    \Sbr{
      \frac{\sqrt{k}|f\Sp k(0)|}{(k-1)!} [\theta
      \rc(f,0)]^{k-1}\times
      n^{-\nth{2k}}\max_{\btw j 1 p} \normx{V_j}_{2k}
    },
  \end{align*}
  and $c_2$ as in \emph{Proposition \ref{prop:lse-d}\/}.
\end{theorem}

If $f$ is linear, then the expression of $c_1$ is simplified into 
\begin{align*}
  c_1 = \sigma\sqrt{2\lambda_p} |f'(0)|
  \max_{\btw j 1 p} \normx{V_j}_2/\sqrt{n}.
\end{align*}
In the general case, as $n^{-1/2k} \max_j\normx{V_j}_k \le \max_j
\normx{V_j}_\infty$,
\begin{align*}
  c_1 \le \sigma\sqrt{2\lambda_p} K
  \max_{\btw j 1 p}\normx{V_j}_\infty, \ \text{with}\
  K=\sum_{k=1}^\infty \frac{\sqrt{k}|f\Sp k(0)|}{(k-1)!} [\theta
  \rc(f,0)]^{k-1}.
\end{align*}
Since $\rc(f,0) = (\Lsup_k |f\Sp k(0)/k!|^{1/k})^{-1}$, it is
easy to see that $c_1<\infty$.  As noted in Section
\ref{ssec:glm-notes}, under mild conditions, $\max_j
\normx{V_j}_\infty = O(\sqrt{\ln (np)})$.  Since $\lambda_p =
O(\ln p)$ and $K$ is a constant, $c_1=O(\sqrt{\ln (np) \ln p})$.
Therefore, for reasonably large $p$, such as $p=n^a$, $c_1=O(\sqrt{\ln
  n})$.  Moreover, as seen previously, under mild conditions, it is
possible that $c_1 = O(\ln n)$.  Combining the comment after
Proposition \ref{prop:lse-d}, it is seen that the regression
estimator \eqref{eq:LSE} can have good precision.

\subsection{Multiple analytic discs}  \label{ssec:multi-disc}
We first need some preparation.  Let $\cN\subset\Coms$ be an open set
containing $I$ such that $f$ has an analytic extension on $\cN$.  Let
$J=\cN\cap\Reals$.  For $u\in \cD(J)$, $i=1,\ldots, n$, and
$k\in\Nats$, define functions,
\begin{align} \label{eq:rc-min}
  a_{i k}(u) = \frac{f\Sp k(X_i\tp u)}{k!},
  \quad
  A_k(u) = \max_{\btw i 1 n} |a_{i k}(u)|,
  \quad
  r(u) = \min_{\btw i 1 n} \rc(f, X_i\tp u).
\end{align}
It is easy to see that $r(u)>0$.  Given any function $b(u)$ on $\cD(J)$
satisfying
\begin{align} \label{eq:radius}
  0<b(u) < r(u)
\end{align}
and given any set $E\subset \cD(J)$, denote
\begin{align} \label{eq:h}
  \sfb(E) = \inf_{u\in E} b(u),
  \quad
  \sfr(E)=\inf_{u\in E} r(u),
  \quad
  \sfA_k(E) = \sup_{u\in E} A_k(u).
\end{align}
If $E$ is finite, then it is easy to see that $\sfr(E) > \sfb(E)$,
and, by $\Lsup_k |a_{i k}(u)|^{1/k} = 1/\rc(f, X_i\tp u)$ for
$u\in\cD(J)$ and $i=1,\ldots, n$,
\begin{align} \label{eq:radius-2}
  \Lsup_{k\toi} \sfA_k(E)^{1/k}
  = \max_{
    \substack{u\in E\\\btw i 1 n}
  } \Lsup_{k\toi} |a_{i k}(u)|^{1/k}
  = \nth{\sfr(E)}.
\end{align}

Let $\grid$ be a subset of $\cD(J)$.  If
\begin{align}   \label{eq:cover}
  E\subset\bigcup_{u\in \grid} \cO_u, \ \text{with}\ \cO_u =
  \ball(u, b(u)/2; \normx{\cdot}_{1,\infty}),
\end{align}
then $\grid$ will be referred to as a ``$b/2$-covering grid'', or
simply ``covering grid'' for $E$.  By this definition, for each point
$u$ in a covering grid and $i=1,\ldots,n$, $f$ is analytic at
$X_i\tp u$ with $\rc(f, X_i\tp u)>b(u)$.  Note that a covering
grid of $E$ need not be its subset. 
If $E$ is compact, it always has a finite covering grid.

Finally, for $E\subset\parsp$, denote
\begin{align*}
  \sfC(E) = \{(1-s) u + sv: s\in [0,1], u, v\in E\},
\end{align*}
i.e., the union of all the line segments connecting pairs of points in
$E$.  If $E$ is bounded (resp.\ compact), then $\sfC(E)$ is bounded
(resp.\ compact).  If $|\sppt(u)|\le a$ for every $u\in E$, then
$|\sppt(v)|\le 2 a$ for every $v\in\sfC(E)$.  However, $\sfC(E)$ may
not be convex, and for unbounded closed $E$, $\sfC(E)$ may not be even
closed.

After all the preparation, the main result can be stated as follows.
\begin{theorem} \label{thm:multi-disc}
  Suppose $I$ is \emph{compact\/} and $\sfd(I,f)>0$.  Fix $\nu\in
  (0,1)$.  In the regression \eqref{eq:LSE}, let $D$ be a closed
  subset of $\pset I {\sfn(\nu)/2}$.  Fix $b(u)$ satisfying
  \eqref{eq:radius}.  Let $\grid$ be a finite $b/2$-covering grid of
  $\sfC(D)$.  Given $q\in (0,1)$, let $\lambda_p=\ln p(1+q^{-1})$.  If
  $\beta\in D$, then the conclusion of \emph{Proposition
    \ref{prop:basic}\/} holds with
  \begin{align} \label{eq:c1-multi-disc}
    c_1 = \sqrt{2}\sigma \sum_{k=1}^\infty
    \Sbr{
      k\sqrt{\ln |\grid| + k\lambda_p}\,
      \sfA_k(\grid) \sfb(\grid)^{k-1}\times n^{-\nth{2k}}
      \max_{\btw j 1 p} \normx{V_j}_{2k}
    }
  \end{align}
  and $c_2$ as in \emph{Proposition \ref{prop:lse-d}\/}.
\end{theorem}

To get $c_1$, it is enough to assume $D$ is a compact subset of
$\cD(J)$.  The stronger assumption that $D\subset \pset
I{\sfn(\nu)/2}$ is needed in order to get both $c_1$ and $c_2$.  By
Proposition \ref{prop:compact}, $\pset I {\sfn(\nu)/2}$ is compact.
Therefore, if $D\subset \pset I{\sfn(\nu)/2}$ is closed, it is compact
as well.

Unlike in Theorem \ref{thm:one-disc}, here $c_1$ depends on $|\grid|$.
In order for the regression estimator \eqref{eq:LSE} to have good
precision, $|\grid|$ has to be controlled.  The smaller $|\grid|$ is,
the higher the precision we can claim for $\est\beta$.  To
see what might be an acceptable level of $|\grid|$, observe that
\begin{align*}
  c_1 \le \sqrt{2}\sigma K \sqrt{\ln|\grid|+\lambda_p}
  \max_{\btw j 1 p}\normx{V_j}_\infty
  = O\Grp{
    \sqrt{\ln(p |\grid|)}\max_j\normx{V_j}_\infty
  },
\end{align*}
where $K = \sum_k k^{3/2}\sfA_k(\grid) \sfb(\grid)^{k-1}$ is finite by
\eqref{eq:radius-2}.  From the comment after Proposition
\ref{prop:lse-d}, it is seen that $\est\beta$ has good precision if
$\sqrt{\ln(p |\grid|)} \max_j\normx{V_j}_\infty = o(\sqrt{n})$.
Provided $\max_j \normx{V_j}_\infty = O(\sqrt{\ln (np)})$ and $p=n^a$,
this implies there should be $\ln |\grid| = o(n/\ln n)$.  Certainly,
$|\grid|$ depends on the choice of the search domain $D$ in
\eqref{eq:LSE} and the property of $f$.  We next get some upper
bounds of $|\grid|$. 

\subsection{Upper bounds on the cardinality of covering grid}
We follow the notation in Section \ref{ssec:multi-disc}.  Recall that
$f$ is analytic on some open domain $\cN\subset\Coms$ containing
$I=[a,b]$ and $J=\cN\cap\Reals$.  The next result says that $|\grid|$
can be as small as 1 in Theorem \ref{thm:multi-disc}.  It follows
directly from the definition of covering grid.
\begin{prop} \label{prop:size-1}
  Let $D\subset B(w, d/2; \, \normx{\cdot}_{1,\infty})$ for some
  $w\in \cD(J)$ and $0<d < r(w)$.   Then for any $b$ satisfying
  \eqref{eq:radius} and $d<b(w)$, $\{w\}$ is a $b/2$-covering grid for
  $\sfC(D)$.
\end{prop}

As an example, if $f$ is analytic in a neighborhood of 0 and 
$\normx{u}_{1, \infty}\le \theta\rc(f,0)/2$ for all $u\in D$, where
$0<\theta<1/2$, then, since $r(0) = \rc(f,0)$, $\{0\}$ is a
$b/2$-covering grid of $\sfC(D)$ for any $b$ satisfying
\eqref{eq:radius} with $b(0)>\theta\rc(f,0)$.

We next consider more general cases.  For ease of notation, for $E
\subset \parsp$ and $S \subset \{1,\ldots,p\}$, denote $\diam(E) =
\diam(E; \normx{\cdot}_{1,\infty})$ and $E_S=\{u\in E: \sppt(u)\subset
S\}$.
\begin{prop} \label{prop:cover-size}
  Fix $b(u)$ satisfying \eqref{eq:radius} and $h\in\Nats$.  Let
  $D\subset \pset I {h/2}$ be compact and $K=\sfC(D)$.

  (1) If $J=\Reals$ and $\bar d_b:=\inf_{u\in \cD(J)} b(u)>0$,
  then $K$ has a $b/2$-covering grid with cardinality no
  greater than
  \begin{align*}
    \sum_{|S|=h:\ K_S\not=\emptyset} \Sbr{2 \diam(K_S)/\bar
      d_b+1}^h 
    \le
    \binom{p}{h} \Sbr{2\diam(D)/\bar d_b + 1}^h.
  \end{align*}

  (2) In general, if $d_b:=\inf_{u\in \cD(I,h)} b(u)>0$, then $K$ has
  a $b/2$-covering grid with cardinality no greater than
  \begin{align*}
    \sum_{|S|=h:\ K_S\not=\emptyset}
    \Sbr{4 \diam(K_S)/d_b+1}^h
    \le
    \binom{p}{h} \Sbr{4\diam(D)/d_b + 1}^h.
  \end{align*}
\end{prop}

Note that, since $I$ is compact, $\inf_{u\in \pset I h} r(u) \ge
\inf_{x\in I} \rc(f, x)>0$, so there are always functions $b(u)$
satisfying \eqref{eq:radius} and $d_b>0$.  For example, $b(u) =
r(u)/2$.

Finally, in Theorem \ref{thm:multi-disc}, $c_1$ depends on the
choice of $\grid$, so it may not be easy to use.   Using the above
bounds on $|\grid|$, we have some more convenient choices for $c_1$,
although they are larger than the one in \eqref{eq:c1-multi-disc}.
\begin{prop} \label{prop:ub}
  Let $D$ be a compact subset of $\pset I{h/2}$ in regression
  \eqref{eq:LSE}.

  (1) Let $\bar d_k = \sup_{x\in J} |f\Sp k(x)|/k!$ and $\bar\rc_0 =
  \inf_{x\in J} \rc(f, x)$.  Suppose $J=\Reals$, $\bar\rc_0>0$, and
  for any $\bar\rc_1\in (0, \bar\rc_0)$, $\sup_{|{\rm Im}(z)|\le
    \bar\rc_1} |f'(z)|<\infty$.  Then the radius of convergence of
  $\sum_{k\ge 1} \bar d_k z^k$ is $\bar \rc_0$ and given $\bar\rc_1\in
  (0,\bar\rc_0)$, $c_1$ in \eqref{eq:c1-multi-disc} can be set equal
  to
  \begin{align}\label{eq:c1-alt1}
    c_1 = \sqrt{2}\sigma \sum_{k=1}^\infty
    \Sbr{
      k\sqrt{ h\ln (p\bar Q)+ k\lambda_p}\,\bar d_k\bar\rc_1^{k-1} 
      \times n^{-\nth{2k}} \max_{\btw j 1 p} \normx{V_j}_{2k}
    },
  \end{align}
  where $\bar Q = 2\diam(D)/\bar\rc_1+1$.

  (2) Let $d_k = \sup_{x\in I} |f\Sp k(x)|/k!$ and $\rc_0=\inf_{x\in
    I} \rc(f, x)$.  Then $\rc_0>0$ is equal to the radius of
  convergence of $\sum_{k\ge 1} d_k z^k$, and given $\rc_1\in
  (0,\rc_0)$, $c_1$ in \eqref{eq:c1-multi-disc} can be set equal to
  \begin{align}\label{eq:c1-alt}
    c_1 = \sqrt{2}\sigma \sum_{k=1}^\infty 
    \Sbr{
      k\sqrt{h\ln(p Q)+k\lambda_p}\,d_k\rc_1^{k-1}
      \times n^{-\nth{2k}} \max_{\btw j 1 p} \normx{V_j}_{2k}
    },
  \end{align}
  where $Q=4\diam(D)/\rc_1+1$.
\end{prop}

In \eqref{eq:c1-alt1}, because the radius of convergence of
$\sum_{k\ge 1} \bar d_k z^k$ is $\bar\rc_0$, $c_1<\infty$.  As
$\lambda_p = \ln[p(1+q^{-1})]$, $c_1 = O(\sqrt{h\ln p} \max_j
\normx{V_j}_\infty)$.  Therefore, under mild conditions, for large
$n$, as long as $h$ is not too large, the regression \eqref{eq:LSE}
still has good precision.

\subsection{Logistic regression with binary noise}  \label{ssec:log-r2}
Let $\eno y n$ be the same random variables as in Section
\ref{ssec:log-r}.  However, we only see their randomly ``flipped''
versions $\eno z n\in \{0,1\}$, such that
\begin{align*}
  \prob\Cbr{\eno z n\gv \eno y n} = \prod_{i=1}^n p_{y_i z_i},
\end{align*}
where $p_{ab}\ge 0$ and $p_{a0}+p_{a1}=1$ for $a=0,1$.  Suppose all
$p_{ab}$ are known.  The regression model now is $\mean(z_i) =
f(X_i\tp \beta)$ with
\begin{align*}
  f(t)=\frac{p_{01} + p_{11} e^t}{1+e^t}.
\end{align*}
If $p_{01}=p_{11}$, then $z_i$ is independent of $y_i$ with
$\prob\Cbr{z_i=1}=p_{11}$, making inference impossible.  Therefore, we
will assume $\Delta_p=|p_{11} - p_{01}|>0$.

Since $f$ is analytic on $\Coms\setminus\{t_k,\ k\in\Ints\}$, where
$t_k = (2k+1)\pi i$, we shall apply Proposition \ref{prop:ub}(1).
First, since $\r_i = z_i - \mean(z_i)$ are independent and $|\r_i|\le
1$, they satisfy the tail assumption \eqref{eq:tail} with $\sigma=1$.
Since $\rc(f, x)$ is  the distance from $z$ to the closest pole,
for any $x\in\Reals$, $\bar\rc_0 = |0-t_1| = \pi$.  Simple
calculation gives $f'(t) = (p_{11} - p_{01}) [2\cosh(t/2)]^{-2}$.
By $2|\cosh(a+bi)| \ge e^{|a|} - e^{-|a|}$ for $a$, $b\in\Reals$, 
it is easy to see that for $y\in (0, \pi)$, 
\begin{align*}
  M(y)
  :=\sup_{|\text{Im} z|\le y} |2\cosh(z/2)|^{-2}<\infty,
\end{align*}

Fix $\bar\rc_1\in (0,\pi)$, $r>0$ and $\nu\in (0,1)$.  Let $I=[-r,r]$ and
\begin{align*}
  D=\Cbr{
    u\in\parsp: \normx{u}_{1,\infty}\le r,\ |\sppt(u)|\le \sfn(\nu)/2
  }.
\end{align*}
Apparently, $D\subset \pset I {\sfn(\nu)/2}$ and $\diam(D) \le r$,
where, as in Proposition \ref{prop:ub}, $\diam(D) =\diam(D;
\normx{\cdot}_{1,\infty})$.

Let $\theta\in (\bar\rc_1/\pi,1)$.  For any $x\in \Reals$ and $k\ge 1$,
by Cauchy's contour integral,
\begin{align*}
  \frac{|f\Sp k(x)|}{k!}
  \le \nth{2k\pi} \oint_{|z-x|=\bar\rc_1/\theta}
  \frac{|f'(z)|\,dz}{(\bar\rc_1/\theta)^k}
  \le
  \frac{\Delta_p\,M(\bar\rc_1/\theta)}
  {k(\bar\rc_1/\theta)^{k-1}},
\end{align*}
giving $\bar d_k \le \Delta_p\,M(\bar\rc_1/\theta) /
[k(\bar\rc_1/\theta)^{k-1}]$.  Therefore, by Proposition
\ref{prop:ub}(1), 
\begin{align*}
  c_1\le\
  &
  \sqrt{2} \Delta_p\, M(\bar\rc_1/\theta)
  \sum_{k=1}^\infty
  \Sbr{\sqrt{R + k\lambda_p}\,
    \theta^{k-1}\times n^{-\nth{2k}} \max_{\btw j 1 p}\normx{V_j}_{2k}
  },
\end{align*}
where $R=\sfn(\nu)/2\times \ln (2 r p/\bar\rc_1+p)$.  On the other
hand, given $r>0$,
\begin{align*}
  m(r):=\inf_{x\in [-r,r]} |2\cosh(x/2)|^{-2}>0.
\end{align*}
Therefore, by Proposition \ref{prop:lse-d},
\begin{align*}
  c_2 \ge \frac{\Delta_p^2 \,m(r)^2 \nu[1+\mu(X)]}{n}
  \min_{\btw j 1 p} \normx{V_j}_2^2.
\end{align*}

Similar to Section \ref{ssec:log-r}, if all the entries of $X$ are
$\pm 1$, then the results can be simplified so that
$D = \Cbr{u\in\parsp: \sum_i |u_i| \le r, \ |\sppt(u)|\le
  \sfn(\nu)/2}$, and
\begin{align*}
  c_1
  \le
  \sqrt{2} \Delta_p\, M(\bar\rc_1/\theta)\sum_{k=1}^\infty
  \sqrt{R + k\lambda_p}\, \theta^{k-1}, \quad
  c_2
  \ge 
  \Delta_p^2\,m(r)^2 \nu[1+\mu(X)].
\end{align*}

\section{Technical details} \label{sec:proofs}
\subsection{Preliminary results}
\begin{proof}[Proof of Proposition \emph{\ref{prop:separable}}]
  (1) Let $S=\sppt(u)$.  If $|S|=0$, then $u=0$ and the inequality
  trivially holds.  Suppose $|S|\ge 1$.  Since $Xu = \sum_{j\in S} u_j
  V_j$,
  \begin{align*}
    \normx{Xu}_2^2 &
    = \sum_{j\in S} |u_j|^2 \normx{V_j}_2^2
    + \sum_{i,\, j\in S, i\not=j} u_i u_j V_i\tp V_j \\
    &
    \ge
    \sum_{j\in S} |u_j|^2 \normx{V_j}_2^2
    - \mu(X) \sum_{i,\, j\in S, i\not=j} |u_i| |u_j|
    \normx{V_i}_2\normx{V_j}_2 \\
    &
    = [1+\mu(X)] \sum_{j\in S} |u_j|^2 \normx{V_j}_2^2
    -\mu(X)\Grp{\sum_{j\in S} |u_j|\normx{V_j}_2}^2.
  \end{align*}
  By Cauchy-Schwartz inequality,
  \begin{align*}
    \normx{Xu}_2^2
    \ge (1+\mu(X) - \mu(X)|S|) \sum_{j\in S} |u_j|^2 \normx{V_j}_2^2.
  \end{align*}
  Since $|S|\le \sfn(\nu) = (1-\nu)[1+1/\mu(X)]$, then $1+ \mu(X)-
  \mu(X)|S| \ge \nu[1+\mu(X)]$, which implies the desired inequality.

  (2) By $\sppt(u-v)\subset \sppt(u)\cup \sppt(v)$ and the assumption,
  $|\sppt(u-v)|\le \sfn(\nu)$.  The inequality then follows from (1).
\end{proof}

\begin{proof}[Proof of Proposition \emph{\ref{prop:compact}}]
  (1) Because $I$ is closed and the mapping $T: u\to Xu$ is
  continuous, $\cD(I) = T^{-1}(I^n)$ is closed.  Also, $\cV_h
  :=\{u\in\parsp: |\sppt(u)|\le h\}$ is closed.  Thus $\pset I h =
  \cD(I) \cap\cV_h$ is closed.  It is easy to see that $\pset I
  h\subset \pset I {h'}$ when $h<h'$.

  (2) Because of (1), to show that $\pset I h$ is compact for
  $h<\sfn(0)$, it suffices to show the set is bounded.  Since
  $h<\sfn(0)$, there is $\nu\in (0,1)$ such that $h\le \sfn(\nu)$.
  Let $u\in \pset I h$.  Then $|\sppt(u)|\le \sfn(\nu)$, so by
  Proposition \ref{prop:separable},
  \begin{align*}
    \normx{u}_2^2 \le \frac{\normx{Xu}_2^2}{\nu(1+\mu(X))
      \min_{\btw j 1 p} \normx{V_j}_2^2
    }.
  \end{align*}
  Since $X_i\tp u\in I$ for each $i$, then $\normx{Xu}_2^2 \le n
  \max_i|X_i\tp u|^2\le n \sup_{x\in I} |x|^2$.  Because
  $I$ is bounded, it is seen $\normx{u}_2^2$ is bounded for  $u\in
  \pset I h$.
\end{proof}

\subsection{Exponential linear models}
In this section, we prove the next two lemma.
\begin{lemma} \label{lemma:glm-c1}
  \emph{Condition \ref{cond:A1}\/} is by satisfied $\varphi=(\eno\varphi
  n)$  with
  \begin{align} \label{eq:c1-glm}
    c_1 = \sigma
    \sqrt{\frac{\ln(p/q)}{2n}} \max_{\btw j 1 p}\normx{V_j}_2.
  \end{align}
\end{lemma}
\begin{lemma} \label{lemma:glm-c2}
  \emph{Condition \ref{cond:A2}\/} is satisfied by $G$ and $\psi =
  (\eno \psi n)$ with
  \begin{align}
    c_2 = 
    \frac{\nu\delta[1+\mu(X)]}{2n} \min_{\btw j 1 p} \normx{V_j}_2^2.
    \label{eq:glm-c2}
  \end{align}
\end{lemma}

By Proposition \ref{prop:basic}, if $c_r = 3c_1^2/c_2$ in
\eqref{eq:MLE}, then \eqref{eq:glm-est} holds with $\kappa_r =
3c_1/c_2$.  Therefore, once the lemmas are proved, we get the
expressions of $c_r$ and $\kappa_r$ as in Theorem \ref{thm:glm}.

As in \eqref{eq:glm-f}, let $G(x) = x_1+\cdots+x_n$ for $x\in\obssp$,
and $\varphi_i(z)=z/2$, $\psi_i(z) = \Lambda(z) -
\Lambda'(X_i\tp\beta) z$ for $\btw i 1 n$ and $z\in\Reals$.

\begin{proof*}[Proof of Lemma \ref{lemma:glm-c1}]
  By \eqref{eq:tail} and $\r\tp V_j = \sum_{i=1}^n X_{ij} \r_i$,
  \begin{align*}
    &
    \prob\Cbr{
      |\r\tp V_j| \le \sqrt{2\ln (p/q)}\,
      \sigma\normx{V_j}_2, \text{ all } j=1,\ldots,p
    } \\
    &
    \ge
    1-
    \sum_{j=1}^p
    \prob\Cbr{
      |\r\tp V_j|^2 > 2\ln (p/q) \sigma^2\normx{V_j}_2^2} 
    \ge
    1-2q.
  \end{align*}
  Consequently, with probability at least $1-2q$,
  \begin{align*}
    \normx{X\tp\r}_\infty = \max_{\btw j 1 p} |\r\tp V_j|
    \le
    \sqrt{2\ln(p/q)}\sigma\max_{\btw j 1 p}\normx{V_j}_2
    = 2 c_1\sqrt{n},
  \end{align*}
  which implies condition \ref{cond:A1} due to the fact that for all
  $u\in\parsp$,
  \begin{align*}
    &
    |\angx{\r}{\varphi(Xu) - \varphi(X\beta)}|
    =
    \nth 2|\angx{\r}{Xu- X\beta}|\\
    &\qquad
    =
    \nth 2|(X\tp \r)\tp(u-\beta)|
    \le \nth 2\normx{X\tp \r}_\infty \normx{u-\beta}_1.
    \tag*{$\Box$}
  \end{align*}
\end{proof*}

\begin{proof}[Proof of Lemma \ref{lemma:glm-c2}]
  Given
  $u\in\pset I {\sfn(\nu)/2}$, for $t\in [0,1]$, let
  \begin{align*}
    h(t) = \sum_{i=1}^n \psi_i((1-t)X_i\tp\beta + t X_i\tp u),
  \end{align*}
  which is well-defined as $(1-t) X_i\tp\beta + t
  X_i\tp u\in I$.  Let $\Delta=G(\psi(Xu) - \psi(X\beta))$.  Then
  \begin{align*}
    \Delta = \sum_{i=1}^n
    \Sbr{\psi_i(X_i\tp u) - \psi_i(X_i\tp\beta)} = h(1)-h(0).
  \end{align*}
  
  Observe that $\psi_i'(X_i\tp \beta) = 0$.  Then $h'(0) = \sum_i
  X_i\tp (u-\beta) \psi_i'(X_i\tp \beta) = 0$, so by Taylor expansion,
  $\Delta= h''(\tau)/2$ for some $\tau\in (0,1)$.  By $\psi_i''(z) =
  \Lambda''(z)$ and $\inf_{t\in I} \Lambda''(t)=\delta>0$,
  \begin{align*}
    \Delta
    &
    = \nth 2\sum_{i=1}^n [X_i\tp (u-\beta)]^2 \psi_i''((1-t) X_i\tp
    \beta+t X_i\tp u) \\
    &
    \ge \frac{\delta}{2}\sum_{i=1}^n [X_i\tp (u-\beta)]^2 = 
    \frac{\delta\normx{X(u-\beta)}_2^2}{2}.
  \end{align*}
  By $|\sppt(u-\beta)|\le |\sppt(u)\cup\sppt(\beta)|\le  \sfn(\nu)$
  and Proposition \ref{prop:separable},
  \begin{align*}
    \Delta
    &\ge \frac{\delta\nu[1+\mu(X)]}{2}
    \sum_{j=1}^p |u_j-\beta_j|^2 \normx{V_j}_2^2
    \ge \frac{\delta\nu[1+\mu(X)]}{2} \min_{\btw j 1 n}
    \normx{V_j}_2^2 \times\normx{u-\beta}_2^2,
  \end{align*}
  and so Condition \ref{cond:A2} is satisfied with $c_2$ set as in
  \eqref{eq:glm-c2}.
\end{proof}

\subsection{Proofs for LS regression: the case of single analytic disc}

First, we establish Condition \ref{cond:A2}.
\begin{proof}[Proof of Proposition \emph{\ref{prop:lse-d}}]
  For $i=1,\ldots, n$ and $u\in D$, since $X_i\tp\beta\in I$ and
  $X_i\tp u\in I$,
  \begin{align*}
    \normx{f(Xu)-f(X\beta)}_2^2
    &= \sum_{i=1}^n |f(X_i\tp u) -  f(X_i\tp\beta)|^2 \\
    &
    \ge
    \sum_{i=1}^n \sfd(f, I)^2 |X_i\tp u-X_i\tp\beta|^2
    = \sfd(f,I)^2 \normx{X(u-\beta)}_2^2.
  \end{align*}
  Since $|\sppt(u-\beta)|\le |\sppt(u)|+|\sppt(\beta)|\le
  \sfn(\nu)$, then by Proposition \ref{prop:separable}, 
  \begin{align*}
    \normx{f(Xu)-f(X\beta)}_2^2
    \ge \sfd(f,I)^2 \nu[1+\mu(X)]\min_{\btw j 1 p}\normx{V_j}_2^2
    \times\normx{u-\beta}_2^2.
  \end{align*}
  Because the right hand side is $c_2 n\normx{u-\beta}_2^2$, the proof
  is complete.
\end{proof}

The main result in this section is Proposition \ref{prop:one-disc},
which together with Proposition \ref{prop:lse-d} immediately leads to
Theorem \ref{thm:one-disc}.  For brevity, in the rest of this section,
we shall denote $\Pi = \{1,\ldots,p\}$.

\subsubsection{Power series expansion  and tail assumption}
To facilitate subsequent discussions, we first consider
\begin{gather*}
  \varphi(x)=(\varphi_1(x_1), \ldots, \varphi_n(x_n))\tp,
  \quad
  x=(\eno x n)\tp\in\obssp,
\end{gather*}
where $\eno\varphi n$ are real-valued functions that may be different
from each other.

Suppose each $\varphi_i$ can be analytically extended to a
neighborhood of $0$ in $\Coms$.  Let
\begin{align}
    a_{i k} = \frac{\varphi_i\Sp k(0)}{k!}\in\Reals,
    \label{eq:power}
\end{align}
Then $\rc(\varphi_i,0) = (\Lsup_k |a_{i k}|^{1/k})^{-1}$.
Since we are interested in $\varphi(Xu) - \varphi(Xv)$ instead of 
$\varphi(Xu)$ itself, without loss of generality, let
$\varphi_i(0)=0$.

For vector $v=(\eno v p)\tp$ and $k$-tuple $\alpha = (\eno\alpha k)
\in\Pi^k$, denote by $v_\alpha$ the product of $v_{\alpha_1}$, \ldots,
$v_{\alpha_k}$.  For example, if $p=3$ and $k=4$, then $v_{(1,3,1,2)}
= v_1 v_3 v_1 v_2= v_1^2 v_2 v_3$.  With this notation, for
$i=1,\ldots, n$, $X_{i\alpha} = X_{i\alpha_1} \cdots X_{i\alpha_k}$.
For each $j=1,\dots, p$, let $n_j(\alpha) = |\{i: \alpha_i=j\}|$.
Clearly, $n_1(\alpha)+\cdots+n_p(\alpha)=k$.

By \eqref{eq:power}, for $i=1,\ldots, n$, provided $|X_i\tp
u|<\rc(\varphi_i,0)$,
\begin{align*}
  \varphi_i(X_i\tp u) = \sum_{k=1}^\infty a_{i k} (X_i\tp u)^k
  = \sum_{k=1}^\infty a_{i k} \Grp{\sum_{\alpha\in \Pi^k} X_{i\alpha}
    u_\alpha}.
\end{align*}
Therefore, if $|X_i\tp u|<\rc(\varphi_i,0)$ for all $i$, then
\begin{align}
  \angx{\r}{\varphi(Xu)}
  = \sum_{i=1}^n \r_i \varphi_i(X_i\tp u)
  = \sum_{k=1}^\infty
  \sum_{\alpha\in \Pi^k}
  \Grp{
    u_\alpha\sum_{i=1}^n \r_i a_{i k} X_{i\alpha}
  }.
  \label{eq:power-series}
\end{align}

\begin{lemma} \label{lemma:control}
  Suppose $\r$ satisfy \eqref{eq:tail}.  Let $q_1$, $q_2$, \ldots $\ge
  0$ with $q:=\sum_k q_k <1/2$.  Given real numbers $\theta_{i k}$,
  $\btw i 1 n$, $\btw k 1 p$, consider the condition
  \begin{align} \label{eq:control}
    \Abs{\sum_{i=1}^n \r_i \theta_{i k} X_{i\alpha}}
    \le \sigma\sqrt{2\ln (p^k/q_k)}
    \sqrt{\sum_{i=1}^n \theta_{i k}^2 X_{i\alpha}^2},
  \end{align}
  where $\sigma$ is the constant in \eqref{eq:tail} and $\ln 0$ is
  defined to be $-\infty$.  Then
  \begin{align}
    &
    \prob\Cbr{
      \text{\eqref{eq:control} holds for \emph{all\/} $k\ge 1$ 
        and $\alpha\in \Pi^k$
      }
    } \ge 1-2q.
    \label{eq:prob-ineq2}
  \end{align}
\end{lemma}
\begin{proof}
  The left hand side of \eqref{eq:prob-ineq2} is at least
  \begin{align*} 
    1-\sum_{k=1}^\infty \sum_{\alpha\in\Pi^k}
    \prob\Cbr{
      \text{\eqref{eq:control} does \emph{not\/} hold for $k$ 
        and $\alpha$
      }
    }
  \end{align*}
  Since $|\Pi^k|=p^k$, it suffices to show that for each $k$ and
  $\alpha=(\eno\alpha k)$,
  \begin{align}  \label{eq:prob-ineq}
    \prob\Cbr{
      \Abs{\sum_{i=1}^n \r_i \theta_{i k} X_{i\alpha}}^2 > 2\sigma^2
      \ln (p^k/q_k) \sum_{i=1}^n \theta_{i k}^2 X_{i\alpha}^2
    } \le 2p^{-k} q_k,
  \end{align}
  which directly follows from \eqref{eq:tail}.
\end{proof}

\subsubsection{Establishing Condition \ref{cond:A1}}
Recall the following multinomial formula: for any $j=1,\ldots, p$,
\begin{align}
  \sum_{\alpha\in\Pi^k} n_j(\alpha) x_j^{n_j(\alpha)-1}
  \prod_{s\not=j} x_s^{n_s(\alpha)}
  = k(x_1+\cdots+x_p)^{k-1},
  \label{eq:mn}
\end{align}
as the left hand side is equal to
\begin{align*}
  &
  \sum_{k_1+\cdots+k_p=k} 
  \binom{k}{k_1\cdots k_p}
  k_j x_j^{k_j-1} \prod_{s\not=j} x_s^{k_s}  \\
  &
  =\frac{\partial}{\partial x_j} \Sbr{
    \sum_{k_1+\cdots+k_p=k} \binom{k}{k_1\cdots k_p} x_j^{k_j}
    \prod_{s\not=j} x_s^{k_s}
  }
  = \frac{\partial\Sbr{(\sum_i x_i)^k}}{\partial x_j}.
\end{align*}

For each $j=1,\ldots, p$, let
\begin{align}  \label{eq:norm-power}
  \omega_{j k} = a_{1 k}^2 X_{1 j}^{2k} + \cdots + a_{n k}^2 X_{n
    j}^{2k}.\\[-2.5em]\nonumber
\end{align}
\begin{lemma} \label{lemma:in-product}
  Suppose that, with $\theta_{i k} = a_{i k}$, \eqref{eq:control}
  holds for \emph{all\/} $k\ge 1$ and $\alpha\in \Pi^k$.   Given $u$
  and $v$, let $d_j = |u_j-v_j|$ and $m_j = |u_j|\vee |v_j|$ for
  $j=1,\ldots, p$.  If
  \begin{align} \label{eq:norm}
    \sum_{j=1}^p m_j \max_{1\le i\le n}
    \frac{|X_{i j}|}{\rc(\varphi_i,0)} < 1,
  \end{align}
  then, letting $\xi=\angx{\r}{\varphi(Xu)-\varphi(Xv)}$,
  \begin{align}
    |\xi|\le
    \sigma\sqrt{2}
    \sum_{k=1}^\infty
    \Sbr{
      k \sqrt{\ln (p^k/q_k)}
      \Grp{
        \sum_{j=1}^p m_j \omega_{j k}^{\nth{2k}}
      }^{k-1} \sum_{j=1}^p d_j \omega_{j k}^{\nth{2k}}
    }.
    \label{eq:in-product-diff2}
  \end{align}
\end{lemma}

\begin{proof}
  By \eqref{eq:norm}, for any $i$,
  \begin{align*}
    |X_i\tp u| \le \sum_{j=1}^p |u_j X_{i j}|
    \le \rc(\varphi_i, 0) \sum_{j=1}^p
    m_j\frac{|X_{ij}|}{\rc(\varphi_i,0)}
    < \rc(\varphi_i,0),
  \end{align*}
  and likewise $|X_i\tp v|<\rc(\varphi_i,0)$.  Therefore, by
  \eqref{eq:power-series},
  \begin{align*}
    \xi = \sum_{k=1}^\infty \sum_{\alpha\in \Pi^k}
    \Sbr{
      (u_\alpha-v_\alpha)
      \sum_{i=1}^n \r_i a_{i k} X_{i\alpha}
    }.
  \end{align*}
  By the assumption, \eqref{eq:control} holds with $\theta_{ik} =
  a_{ik}$ for all $k\ge 1$ and $\alpha\in\Pi^k$.  Thus
  \begin{align}
    |\xi|
    &
    \le
    \sum_{k=1}^\infty \sum_{\alpha\in \Pi^k}
    \Sbr{
      |u_\alpha-v_\alpha|
      \Abs{\sum_{i=1}^n \r_i a_{i k} X_{i\alpha}}
    } \nonumber\\
    &
    \le
    \sum_{k=1}^\infty 
    \sigma\sqrt{2\ln(p^k/q_k)}
    \Sbr{
      \sum_{\alpha\in \Pi^k}|u_\alpha-v_\alpha|\sqrt{M_\alpha}
    }
  \label{eq:in-product-diff}
  \end{align}
  where $M_\alpha = \sum_{i=1}^n a_{i k}^2 X_{i\alpha}^2$. 
  Given $k\ge 1$, for each $\alpha\in\Pi^k$, by
  $n_1(\alpha)+\cdots+n_p(\alpha)=k$ and Cauchy-Schwartz inequality,
  \begin{align*}
    M_\alpha=
    \sum_{i=1}^n a_{i k}^2 \prod_{j=1}^p X_{i j}^{2n_j(\alpha)} 
    \le \prod_{j=1}^p \Grp{\sum_{i=1}^n a_{i k}^2
      X_{i j}^{2k}}^{n_j(\alpha)/k}
    \le \prod_{j=1}^p \omega_{j k}^{n_j(\alpha)/k},
  \end{align*}
  where the last inequality is due to the notation in
  \eqref{eq:norm-power}.  On the other hand,
  \begin{align*}
    |u_\alpha - v_\alpha| &
    = \Abs{
      \prod_{j=1}^p u_j^{n_j(\alpha)} - \prod_{j=1}^p v_j^{n_j(\alpha)}
    } \\
    &
    \le \sum_{j=1}^p \Cbr{
      |u_j^{n_j(\alpha)} - v_j^{n_j(\alpha)}|
      \prod_{s=1}^{j-1} |v_s|^{n_s(\alpha)} \prod_{s=j+1}^n
      |u_s|^{n_s(\alpha)}
    } \\
    &
    \le
    \sum_{j=1}^p n_j(\alpha) d_j m_j^{n_j(\alpha)-1} 
    \prod_{s\not=j} m_s^{n_s(\alpha)}.
  \end{align*}
  Therefore,
  \begin{align*}
    &\hspace{-2ex}
    \sum_{\alpha\in\Pi^k}
    |u_\alpha-v_\alpha| \sqrt{M_\alpha}
    \le
    \sum_{\alpha\in\Pi^k}
    \Cbr{
      \Sbr{
        \sum_{j=1}^p 
        n_j(\alpha) d_j m_j^{n_j(\alpha)-1}
        \prod_{s\not=j} m_s^{n_s(\alpha)}
      }
      \prod_{j=1}^p \omega_{j k}^{n_j(\alpha)/(2k)}
    }
    \\
    &\qquad
    =
    \sum_{j=1}^p d_j \omega_{j k}^{\nth{2k}}
    \Cbr{
      \sum_{\alpha\in\Pi^k} n_j(\alpha)
      \Sbr{m_j\omega_{j k}^{\nth{2k}}}^{n_j(\alpha)-1}
      \prod_{s\not=j} \Sbr{m_s \omega_{s k}^{\nth{2k}}}^{n_s(\alpha)}
    }\\
    &\qquad
    = k\Grp{
      \sum_{j=1}^p m_j \omega_{j k}^{\nth{2k}}
    }^{k-1} \sum_{j=1}^p d_j \omega_{j k}^{\nth{2k}},
  \end{align*}
  where the last equality is due to the multinomial formula
  \eqref{eq:mn}.  Now by \eqref{eq:in-product-diff}, the inequality in
  \eqref{eq:in-product-diff2} is proved.
\end{proof}

\begin{prop} \label{prop:one-disc}
  Fix $\theta \in (0,1)$.  Let $D=\Cbr{u\in\parsp:
    \normx{u}_{1,\infty}\le \theta\rc(f,0)/2}$ in \emph{Condition
    \ref{cond:A1}\/} and $\r$ satisfy \eqref{eq:tail} for
  $\sigma>0$.   If $\beta\in D$, then \emph{Condition \ref{cond:A1}\/}
  is satisfied by setting $c_1$ as in \emph{Theorem
    \ref{thm:one-disc}\/}. 
\end{prop}

\begin{proof}
  We have $\varphi_i=f$ and $\rc(\varphi_i,0)=\rc(f,0)$.  For $u\in
  D$, let $d=u-\beta$ and $m=(\eno m n)\tp$, with $m_j=|u_j|\vee
  |\beta_j|$.  Then
  \begin{align} \label{eq:w-L1}
    \normx{m}_{1,\infty} \le \normx{u}_{1,\infty} +
    \normx{\beta}_{1,\infty} \le \theta\rc(f,0).
  \end{align}
  As a result
  \begin{align*} 
    \sum_{j=1}^p m_j \max_{\btw i 1 n} \frac{|X_{ij}|}{\rc(f,0)}
    = \frac{\normx{m}_{1,\infty}}{\rc(f,0)} \le \theta
  \end{align*}
  and \eqref{eq:norm} is satisfied.  Let $q_k = (\frac{q}{1+q})^k$.
  Then $\sum_k q_k=q$, so by Lemmas \ref{lemma:control} and
  \ref{lemma:in-product}, with probability at least $1-2q$,
  \eqref{eq:in-product-diff2} holds.  For each $k\ge 1$, by the
  notation in \eqref{eq:norm-power}, $\omega_{j k} = (|f\Sp k
  (0)|/k!)^2 \normx{V_j}_{2k}^{2k}$.  Recall that in Theorem
  \ref{thm:one-disc}, $\lambda_p$ is defined to be $\ln[p(1+q^{-1})]$.
  Since $\sqrt{\ln (p^k/q_k)} = \sqrt{k\lambda_p}$,
  \begin{align} \label{eq:w-MD}
    k\sqrt{\ln (p^k/q_k)}\Grp{
      \sum_{j=1}^p m_j \omega_{j k}^{\nth{2k}}\!
    }^{k-1}\!\! \sum_{j=1}^p d_j \omega_{j k}^{\nth{2k}}
    =
    \frac{\sqrt{k\lambda_p}|f\Sp k(0)|}{(k-1)!}
    \times\normx{m}_{1,2k}^{k-1}
    \normx{d}_{1,2k},
  \end{align}
  where the weighted $L_1$ norm $\normx{\cdot}_{1,s}$ is defined in
  \eqref{eq:w-L1-norm} and satisfies
  \begin{align*}
    \normx{u}_{1,s}\le
    \begin{cases}
      n^{1/s} \normx{u}_{1,\infty}, \\
      \displaystyle
      \max_{\btw j 1 p} \normx{V_j}_s\times \normx{u}_1,
    \end{cases}\quad
    s\ge 1.
  \end{align*}
  Then by \eqref{eq:w-L1},
  \begin{align*}
    \normx{m}_{1,2k}^{k-1}\normx{d}_{1,2k}
    &
    \le \Grp{n^{\nth{2k}} \normx{m}_{1,\infty}}^{k-1}\times
    \max_{\btw j 1 p}\normx{V_j}_{2k} \times \normx{d}_1\\
    &
    \le \sqrt{n}\,[\theta\rc(f,0)]^{k-1} \times
    n^{-\nth{2k}}
    \max_{\btw j 1 p} \normx{V_j}_{2k}\times
    \normx{d}_1.
  \end{align*}
  Together with \eqref{eq:in-product-diff2} and \eqref{eq:w-MD}, this
  yields the proof.
\end{proof}

\subsection{LS regression: multiple analytic disc  case}
\subsubsection{Proof of  Theorem \ref{thm:multi-disc}}
We first restate Lemma \ref{lemma:control} as follows.
\begin{lemma} \label{lemma:control-m}
  Let $\r$ satisfy \eqref{eq:tail}.  Let $E\subset \cD(J)$ be finite
  and for $k\ge 1$ and $u\in E$, let $q_{k,u}\ge 0$, such that
  $q:=\sum_k \sum_{u\in E} q_{k,u} <1/2$.  Consider the condition
  \begin{align} \label{eq:control-m}
    \Abs{\sum_{i=1}^n \r_i a_{i k}(u) X_{i\alpha}}
    \le \sigma\sqrt{2\ln (p^k/q_{k,u})}
    \sqrt{\sum_{i=1}^n a_{i k}(u)^2 X_{i\alpha}^2},
  \end{align}
  where $\sigma>0$ is the constant in \eqref{eq:tail}.  Then
  \begin{align*}
    &
    \prob\Cbr{
      \text{\eqref{eq:control-m} holds for \emph{all\/} $k\ge 1$,
        $\alpha\in \Pi^k$, and $u\in E$
      }
    } \ge 1-2q.
  \end{align*}
\end{lemma}

The next result provides a bound on $|\angx{\r}{f(Xu)-f(Xv)}|$ for
suitable $u$ and $v$.  The method of its proof is describe at the end
of Section \ref{sec:LSE-prelim}.
\begin{lemma} \label{lemma:ip-m}
  Given  $b(u)$ satisfying \eqref{eq:radius}, let $\grid$ be a finite
  $b/2$-covering grid of a set $K\subset \cD(J)$.  Fix $q_k\ge 0$ such
  that $q:=\sum_k q_k <1/2$ and $\ln q_k = O(k)$ over
  $\cJ=\{k\in\Nats:  \sfA_k(\grid)>0\}$.  Suppose that, with $E=\grid$
  and $q_{k,u}= q_k/|G|$, \eqref{eq:control-m} holds for all $k\ge 1$,
  $\alpha\in \Pi^k$, and $u\in \grid$.  If $u$, $v\in K$ and the
  entire line segment connecting them is in $K$, then, letting $\xi =
  \angx{\r}{f(Xu)-f(Xv)}$ and $d=v-u$,
  \begin{align} \label{eq:ip-m2}
    |\xi|\le \sigma \sqrt{2n}\,H(\sfb(\grid), d)
  \end{align}
  where $H(\sfb(\grid), d) < \infty$, with
  \begin{align*}
    H(z, d)=
    \sum_{k=1}^\infty
    k\sqrt{\ln |G|+\ln(p^k/q_k)}\ \sfA_k(\grid) \times
    n^{-\nth{2k}} \normx{d}_{1,2k}\times
    z^{k-1}.
  \end{align*}
\end{lemma}

\begin{proof}
  Since $\grid$ is finite, $\sfb(\grid)<\sfr(\grid)$.  Given $\rr\in
  (0, \sfr(\grid)/\sfb(\grid)-1)$, let
  \begin{align*}
    T=\left\lceil{
        \frac{2\normx{d}_{1,\infty}}{\rr \sfb(\grid)}
      }\right\rceil.
  \end{align*}
  By the assumption, $u + \theta d\in K$ for $\theta\in [0,1]$.
  For $t=0, \ldots, T$, let $u\Sp t = u+t d/T$.  Then $u\Sp 0 =
  u$, $u\Sp T=v$, and $u\Sp t\in K$.  Fix $t=1$, \ldots,
  $T$.  Then
  \begin{align*}
    \begin{split}
      \normx{u\Sp t - u\Sp{t-1}}_{1,\infty}
      = \normx{d}_{1,\infty}/T
      \le \rr \sfb(\grid)/2.
    \end{split}
  \end{align*}

  By the definition of $\grid$, we can find some $w\in \grid$, such
  that $\normx{u\Sp t - w}_{1,\infty} \le \sfb(\grid)/2$.  Then
  $\normx{u\Sp {t-1} - w}_{1,\infty} \le (1+\rr) \sfb(\grid)/2$.  Let
  $\varphi(x) = (\varphi_1(x_1), \ldots, \varphi_n(x_n))\tp$, with
  \begin{align*}
    \varphi_i(z) = f(z+X_i\tp w) - f(X_i\tp w), \quad\btw i 1 n.
  \end{align*}
  Let $\tilde u = u\Sp t-w$, $\tilde v=u\Sp{t-1}-w$.  Then
  \begin{align*}
    \varphi(X\tilde u) = f(Xu\Sp t) - f(Xw), 
    \quad
    \varphi(X\tilde v) = f(Xu\Sp {t-1}) - f(Xw),
  \end{align*}
  and, as shown just now,
  \begin{align*} 
    \normx{\tilde u}_{1,\infty}\le (1+\rr) \sfb(\grid)/2,
    \quad
    \normx{\tilde v}_{1,\infty}\le (1+\rr) \sfb(\grid)/2.
  \end{align*}

  Let $m=(\eno m p)\tp$ with $m_j = |\tilde u_j|\vee |\tilde v_j|$.
  From the above equalities we get
  \begin{align} \label{eq:m-1-infty}
    \normx{m}_{1,\infty} \le 
    \normx{\tilde u}_{1,\infty} + \normx{\tilde v}_{1,\infty} 
    \le (1+\rr)\sfb(\grid),
  \end{align}
  and hence, by $\rc(\varphi_i, 0) = \rc(f, X_i\tp w)
  \ge r(w)$,
  \begin{align*}
    \sum_{j=1}^p m_j \max_{\btw i 1 n}
    \frac{|X_{i j}|}{\rc(\varphi_i,0)}
    \le\frac{\normx{m}_{1,\infty}}{r(w)}
    \le
    \frac{(1+\rr)\sfb(\grid)}{\sfr(\grid)}<1.
  \end{align*}

  Now Lemma \eqref{lemma:in-product} can be applied to $\varphi$,
  with $u$, $v$, and $q_k$ therein replaced with $\tilde u$, $\tilde
  v$, and $q_k/|\grid|$, respectively.  Then
  \begin{align*}
    \Abs{\angx{\r}{f(Xu\Sp t)-f(Xu\Sp{t-1})}} 
    \le
    \sigma\sqrt{2}
    \sum_{k=1}^\infty k \sqrt{\ln(|\grid|p^k/q_k)}\,M_k\Sp t,
  \end{align*}
  where
  \begin{align*}
    \begin{split}
      M_k\Sp t = \Grp{
        \sum_{j=1}^p m_j \omega_{j k}^{\nth{2k}}
      }^{k-1}
      \sum_{j=1}^p |u\Sp t_j - u\Sp{t-1}_j|
      \omega_{j k}^{\nth{2k}},
      \hspace{1cm}\\
      \text{with}\quad
      \omega_{j k}
      =\sum_{i=1}^n a_{i k}^2(w)|X_{i j}|^{2k}
      \le \sfA_k^2(\grid) \normx{V_j}_{2k}^{2k}
      \le n \sfA_k^2(\grid)\normx{V_j}_\infty^{2k}.
    \end{split}
  \end{align*}
  Since $u\Sp t-u\Sp {t-1} = d/T$, it follows that
  \begin{align*}
    M_k\Sp t
    &\le \Grp{
      \sum_{j=1}^p
      m_j n^\nth{2k} \sfA_k^\nth{k}(\grid)\normx{V_j}_\infty
    }^{k-1} \times \frac{\sfA_k^\nth k(\grid)}{T}
    \sum_{j=1}^p d_j \normx{V_j}_{2k} \\
    &
    = \frac{\sqrt{n} \sfA_k(\grid)}{T} \normx{m}_{1,\infty}^{k-1} 
    \times n^{-\nth{2k}}\normx{d}_{1,2k} \\
    &
    \le 
    \frac{\sqrt{n} \sfA_k(\grid)}{T} [(1+\rr)\sfb(\grid)]^{k-1}
    \times n^{-\nth{2k}} \normx{d}_{1,2k},
  \end{align*}
  where the last inequality is due to \eqref{eq:m-1-infty}.
  Consequently,
  \begin{align*}
    |\xi|
    &\le
    \sum_{t=1}^T \Abs{\angx{\r}{f(Xu\Sp t)-f(Xu\Sp{t-1})}} 
    = H((1+\rr) \sfb(\grid), d).
  \end{align*}

  By \eqref{eq:radius-2} and $\ln q_k=O(k)$ over $\cJ$, the radius
  of convergence of the power series defining $g(z)=H(z,d)$ is
  $\sfr(\grid) > \sfb(\grid)$.  As $(1+\rr)\sfb(\grid)<\sfr(\grid)$,
  we can let $\rr\to 0$ and apply dominated convergence.  The proof is
  then complete.
\end{proof}

\begin{prop} \label{prop:multi-disc}
  In \emph{Condition \ref{cond:A1}\/}, let $D$ be a compact subset of
  $\cD(J)$.  Suppose $\r$ satisfies \eqref{eq:tail} for some
  $\sigma>0$.  Let $\grid$ be a finite $b/2$-covering grid of
  $\sfC(D)$.  If $\beta\in D$, then \emph{Condition \ref{cond:A1}\/}
  is satisfied by setting $c_1$ as in \emph{Theorem
    \ref{thm:multi-disc}\/}.
\end{prop}
\begin{proof}
  Since $\sfC(D)$ is compact, it indeed has a finite $b/2$-covering
  grid, justifying the assumption on $\grid$.  As in the proof of
  Proposition \ref{prop:one-disc}, let $q_k = (\frac{q}{1+q})^k$.
  Then by Lemmas \ref{lemma:control-m} and \ref{lemma:ip-m}, with
  probability at least $1-2q$, \eqref{eq:ip-m2} holds.  The rest of
  the proof follows that for Proposition \ref{prop:one-disc} and hence
  is omitted for brevity.
\end{proof}

\begin{proof}[Proof of Theorem \emph{\ref{thm:multi-disc}}]
  First, by $D\subset \pset I {\sfn(\nu)/2}$ and $\sfd(I, f)>0$, 
  Proposition \ref{prop:lse-d} can be applied to yield $c_2$.  Second,
  $\sfC(D)$ is compact and since $I$ is an interval, $\sfC(D)\subset
  \cD(I)$.  Then $\sfC(D) \subset \cD(J)$.  Proposition
  \ref{prop:multi-disc} can be applied to $K=\sfC(D)$ to get $c_1$.
\end{proof}

\subsubsection{Other technical results}
\begin{proof}[Proof of Proposition \emph{\ref{prop:cover-size}}]
  Because $D\subset \pset I {h/2}$ and is compact, $K=\sfC(D)\subset
  \pset I h$ and is compact.

  First, fix $S$ with $|S|=h$ and $K_S\not= \emptyset$.  Let
  $\psi_S: \parsp\to \Reals^S$ be the natural projection and 
  $\imath_S: \Reals^S\to \parsp$ the immersion, such that
  $\imath_S(y) = z\in\parsp$, with $z_j=y_j$ for $j\in S$ and $z_j=0$
  for $j\not\in S$.  Define the weighted $L_1$ norm $\normx{\cdot}_S$
  on $\Reals^S$ such that $\normx{u}_S= \sum_{j\in S}
  |u_j|\normx{V_j}_\infty$.  For ease of notation, denote
  $B_S(w,a)=B(w,a; \normx{\cdot}_S)$ and $\delta_S(E) = \delta(E;
  \normx{\cdot}_S)$.  Likewise, denote $B(w,a) = B(w,a;
  \normx{\cdot}_{1,\infty})$ and $\delta(E) = \delta(E;
  \normx{\cdot}_{1,\infty})$. 

  Fix $d>0$.  Later we will set $d$ to specific values.  Let
  $E=\psi_S(K_S)$.  It is easy to verify that $\delta_S(E) =
  \delta(K_S)$.  By simple geometric argument, it is seen that $E$ can
  be covered by no more than $\Sbr{\diam(K_S)/d+1}^h$ spheres
  $B_S(\tilde u_k, d)$, with each one intersecting with $E$.  Let $u_k
  = \imath_S(\tilde u_k)$.

  In case (1), let $d=\bar d_b/2$.  By $J=\Reals$, $f$ is
  analytic at every $X_i\tp u_k$.  Then, by
  \begin{align*}
    K_S = \imath_S(E)
    \subset \bigcup_k \imath_S(B(\tilde u_k, d))
    \subset \bigcup_k B(u_k,d)
    \subset \bigcup_k B(u_k,b(u_k)/2),
  \end{align*}
  $\eno u m$ is a $b/2$-covering grid of $K_S$.

  In case (2), Let $d=d_b/4$.  Since $f$ may not be analytic at
  every $X_i\tp u_k$, we cannot directly take $\eno u m$ as a covering
  grid.  For each $i=1,\ldots,m$, choose an arbitrary $\tilde w_k \in
  B_S(\tilde u_k, d)\cap E$ and let $w_k = 
  \imath_S(\tilde w_k)$.  As $w_k \subset K_S$, $f$ is
  analytic at every $X_i\tp w_k$.  It is easy to check that
  $B_S(\tilde w_k, 2d)$ contains $B_S(\tilde u_k, d)$.  Therefore,
  \begin{align*}
    K_S = \imath_S(E)
    \subset \bigcup_k \imath_S(B(\tilde w_k, 2d))
    \subset \bigcup_k B(w_k, 2d)
    \subset \bigcup_k B(w_k,b(u_k)/2),
  \end{align*}
  so $\eno w m$ is a $b/2$-covering grid of $K_S$.

  Denote by $\grid_S$ the covering grid as above in either
  case.  As $K=\bigcup_{|S|=h} K_S$, $\grid = \bigcup_{|S|=h:
    K_S\not=\emptyset} G_S$ is a $b/2$-covering grid of $K$ and
  \begin{align*}
    |\grid| \le \sum_{|S|=h: K_S\not=\emptyset} |\grid_S| 
  \end{align*}
  We already know $|\grid_S| \le [\diam(K_S)/d+1]^h$.  By
  $\diam(K_S)\le \diam(K) = \diam(D)$, 
  \begin{align*}
    |\grid_S|\le [\diam(D)/d+1]^h.
  \end{align*}
  Finally, there are at most $\binom{p}{h}$ subsets $S$ with $|S|=h$
  and $K_S \not=\emptyset$.  The proof for the bounds on $|\grid|$ is
  thus complete.
\end{proof}

\begin{proof}[Proof of Proposition {\rm \ref{prop:ub}}]
  (1) If $c>\bar\rc_0$, then there is $t\in J$ such that $\rc(f,t)<c$.
  Since $\Lsup_k |f\Sp k(t)/k!|^{1/k} = \rc(f,t)^{-1}$,
  \begin{align*}
    \Lsup_{k\toi} \bar d_k c^k
    \ge \Lsup_{k\toi} \frac{|f\Sp k(t)| c^k}{k!}
    =\infty.
  \end{align*}
  Therefore, the radius of convergence of $\sum_{k\ge 1} \bar d_k z^k$
  is at most $\bar\rc_0$.  To show that the radius of convergence is
  $\bar\rc_0$, it suffices to show that $\bar d_k c^k$ is bounded for
  any $c\in (0,\bar\rc_0)$.  By assumption $M := \sup_{|{\rm
      Im}(z)|\le c}|f'(z)|<\infty$.  Fix $x\in \Reals$.  For any $z$
  with $|z-x|=c$, $|{\rm Im}(z)|\le c$.  Therefore, by Cauchy's
  contour integral,
  \begin{align*}
    \frac{|f\Sp k(x)|}{k!}
    &= \Abs{\frac{1}{2k\pi\sqrt{-1}} \oint_{|z-x|=c}
      \frac{f'(z) d z}{(z-x)^k}
    }
    \le \frac{1}{2k\pi} \oint_{|z-x|=c}
    \frac{|f(z)| dz}{|z-x|^k} \le \frac{M}{k c^{k-1}}.
  \end{align*}
  Take supremum over $x\in\Reals$.  Then we get $\bar d_k c^k \le
  M/c<\infty$ for all $k\ge 1$.

  From the definitions in \eqref{eq:rc-min}, it is clear that $A_k(u)
  \le \bar d_k$ and $r(u) \ge \bar \rc_0$ for $u\in D$.  Given
  any $\bar\rc_1\in (0, \bar\rc_0)$, let $b(u)\equiv \bar\rc_1$.  
  By Proposition \ref{prop:cover-size}~(1), there is a $b/2$-covering
  grid $\grid$ for $\sfC(D)$ with $|\grid|\le p^h
  (2\diam(D)/\bar\rc_0+1)^h$.  Therefore, $c_1$ can be set as in
  \eqref{eq:c1-alt1}.

  (2)  For each $x\in I$, $\rc(f,x)>0$.  Since $I$ is compact,
  it is covered by a finite number of intervals $(x_i - \rc(f,x_i)/2,
  x_i + \rc(f,x_i)/2)$.  Let $c=\min_i \rc(f,x_i)/2$.  Then $c>0$.
  For any $x\in I$, there is $x_i$ such that $|x-x_i|< \rc(f,x_i)/2$.
  Then for any $z\in \Coms$ with $|z-x|<c$, $|z-x_i|< \rc(f,x_i)$ and
  hence $f$ is analytic at $z$.  As a result, $f$ is analytic in the
  disc centered at $x$ with radius $c$, and so $\rc(f, x)\ge c$.
  This leads to $\rc_0=\inf_{x\in I} \rc(f,x)\ge c$.  For $c\in (0,
  \rc_0)$, since $I_c = \{z\in\Coms: |z-x|\le c$ for some
  $x\in I\}$ is compact, $M:=\sup_{z\in I_c}|f'(z)|<\infty$.
  Using Cauchy's contour integral as in (1), it can be shown that
  $\rc_0$ is the radius of convergence of $\sum_{k\ge 1} d_k z^k$.
  The rest of (2) can be proved following the argument for (1).
\end{proof}
  
\bibliographystyle{agsm}
\bibliography{ldpdb,Estimate,NSdb.bib}

\end{document}